\documentclass[11pt]{amsart}
\usepackage{epsfig,amssymb,rotating}

\newtheorem{theorem}{Theorem}[section]
\newtheorem{lemma}[theorem]{Lemma}
\newtheorem{corollary}[theorem]{Corollary}
\newtheorem{proposition}[theorem]{Proposition}

\newtheorem{problem}[theorem]{Problem}

\newtheorem{definition}[theorem]{Definition}
\newtheorem{example}[theorem]{Example}

\newtheorem{question}[theorem]{Question}

\newcommand{\diag}[2]{\parbox{#2}{\psfig{figure=#1.eps,height=#2}}}
\newcommand{\comment}[1]{\,}
\newcommand{\cal}{\mathcal}
\newcommand{\p}{\partial}
\newcommand{\A}{\mathbb A}
\newcommand{\R}{\mathcal R}

\newcommand{\E}{\mathcal E}
\newcommand{\V}{\mathcal V}

\newcommand{\Rb}{\mathbb R}
\newcommand{\C}{\mathbb C}
\renewcommand{\S}{\mathcal S}
\newcommand{\G}{\mathcal G}
\newcommand{\W}{\mathcal W}
\newcommand{\Z}{\mathbb Z}
\newcommand{\B}{\mathbb B}

\newcommand{\ve}{\varepsilon}
\newcommand{\ld}{{\cal L}{\cal D}}

\title{Confluence Theory for Graphs}

\author{Adam S. Sikora, Bruce W. Westbury}

\thanks{The first author was sponsored in part by NSF grant
\#DMS-0307078.}

\begin{document}

\thispagestyle{empty}

\begin{abstract} 
We develop a theory of confluence of
graphs. We describe an algorithm for proving that a given system of
reduction rules for abstract graphs and graphs in surfaces is locally 
confluent.
We apply this algorithm to show that each simple Lie algebra of rank at most 
$2$, gives rise to a confluent system of reduction rules of 
graphs (via Kuperberg's spiders) in an arbitrary surface. 
As a further consequence of this result, we find
canonical bases of $SU_3$-skein modules of cylinders over orientable surfaces.
\end{abstract}

\pagestyle{myheadings} 

\maketitle
\tableofcontents \vspace{.2in}

%
\section{Introduction}
\label{s_intro}
%

This paper is motivated by the following problem appearing in
representation theory of Lie algebras and of quantum groups, in the study
of moduli spaces, in knot theory, and in other areas of mathematics.
We state it first for abstract graphs and, later, for graphs in manifolds.

Let $R$ be a ring. An {\em $R$-linear graph} is a formal
$R$-linear combination of graphs 
$\Gamma=\sum_{i=1}^k r_i \Gamma_i,$ such that the graphs
$\Gamma_i$ have distinguished sets $E_i$ of $1$-valent vertices
(called {\em external}) and there are specified bijections 
$E_1\simeq E_2\simeq .... \simeq E_k$.
For any graph $\Gamma'$ with a distinguished set of 
$1$-valent external vertices $E'$ in a bijection with $E_1$ 
(and, consequently, in a bijection with $E_i$ for all $i$), 
let $<\Gamma_i,\Gamma'>$ denote the contraction of $\Gamma_i$ and 
$\Gamma'$ along their external vertices, respecting the specified 
bijections.
In the process of the contraction these $1$-valent vertices are removed 
and adjacent edges identified.
Finally, let $<\Gamma,\Gamma'>= \sum_{i=1}^k r_i <\Gamma_i,\Gamma'>$.

Let $\G$ be a set of graphs, $\{\Gamma_i\}_{i\in I}$ be a set of 
$R$-linear graphs, and let $\R(\Gamma_i,i\in I)\subset R\G$ be the submodule 
generated by contractions $<\Gamma_i,\Gamma'>$ for all $i\in I$ and 
all graphs $\Gamma'$ as above. 
\begin{enumerate}
\item Is $R\G/{\mathcal R}(\Gamma_i, i\in I)$ a free $R$-module?
If so, then find an explicit basis of it. 
\item Can a basis be given by taking all graphs in $\G$ satisfying 
a certain ``natural'' property?
\end{enumerate}
Examples appear in Section \ref{s_dichrom}.

The topological version of this problem in dimension $n$ 
involves topological graphs embedded in $n$-dimensional manifolds.
An {\em $R$-linear topological graph} is $\Gamma=\sum_i^k r_i\Gamma_i,$ 
such that $\Gamma_1,...,\Gamma_k$ lie in a manifold $M$ of dimension $n$ 
and there is a finite set
$E\subset \p M,$ such that $\Gamma_i\cap \p M=E$ for every $i$ 
and this set is composed of 
$1$-valent vertices of $\Gamma_i$. 
If $\imath: M\to N$ is an embedding into a manifold of equal dimension and
$\Gamma'$ is a graph in $\overline{N\setminus \imath(M)}$
such that points of $E$ are $1$-valent vertices of $\Gamma'$
then $<\Gamma_i,\Gamma'>$ denotes the contraction 
of graphs $\Gamma_i$ and $\Gamma'$ along the vertices in $E$.
As before, $<\Gamma,\Gamma'>=\sum_i r_i<\Gamma_i,\Gamma'>$.

Now, let $\G$ be a set of topological graphs in $N$ and
$\Gamma_i$ be an $R$-linear graph in $M_i,$ for every $i$ in some index 
set $I$.
As before, let $\R(\Gamma_i,i\in I)\subset R\G$ be the submodule generated by 
$<\Gamma_i,\Gamma'>$ for all $i$'s and all embeddings $\imath:M_i\to N$
and all graphs $\Gamma'$ as above.
In this setting we ask again questions (1),(2) above. 

The flavor of these questions depends on the 
dimension of the manifold $N$:

\noindent{\bf (Dim=2)} 
Interesting examples come from Kuperberg's spider webs, \cite{Ku-spider}, 
which provide a convenient graphical description of representation 
theory of Lie algebras and associated quantum groups of rank $\leq 2.$
These are spaces of graphs in $D^2$ considered modulo certain 
relations, of the type defined above.
The classes of graphs considered and the relations between them
depend on the Lie algebra in question.
Because of their applications to quantum invariants, it is important to
consider Kuperberg's webs in surfaces other than 
$D^2$ as well. We answer (1),(2) for these graphs in 
Sections \ref{s_a1}-\ref{s_g2}. 
Our approach is based on theory of confluence of graphs developed in 
Section \ref{s_confluence} and an algorithm for finding confluent 
reduction rules for graphs described in Section \ref{ss_overlaps}.
As an application, we will find canonical bases of skein modules of 
skein modules of $[0,1]$-bundles over surfaces for all simple Lie 
groups of rank $1$ and $2.$ This reproves theorem of 
Przytycki, \cite[Thm 3.1]{P-fundamentals},
for the Kauffman bracket ($SU_2$) skein modules of 
$[0,1]$-bundles over surfaces and answers the question 
for $SU_3$-skein modules, c.f. \cite{FZ,S-SUn}.

\noindent{\bf (Dim=3)} The three-dimensional version of this problem 
appears in knot theory, for example, in connection with 
Vassiliev invariants and skein modules. In both cases,
(1),(2) are open in general.

\noindent{\bf (Dim$>$3)}
In dimensions greater than $3$ homotopic graphs are isotopic,
and therefore the problem of describing $R\G/{\mathcal R}(\Gamma_i, i\in I)$
can be reduced to purely algebraic form depending on $\pi_1(N)$ only, since 
every $\Gamma\subset N$ is determined by a labeling all cycles of $\Gamma$ by
conjugacy classes of $\pi_1(N)$.
In particular, if $\pi_1(N)$ is trivial then 
graphs in $N$ can be thought as abstract graphs.
For that reason, it is enough to consider questions (1),(2) for abstract
graphs only.

%
\section{Confluence}\label{s_confluence}
%

We will approach the problems outlined in Introduction, by the method of 
confluence. To introduce it in its most abstract form, consider a set 
of objects $V$ and a set of reduction rules, $E$, composed of pairs 
of elements of $V,$ denoted by $v\to v'$. In other words, $(V,E)$ is 
an arbitrary directed graph. A sequence of its vertices 
$v_1\to v_2\to ... \to v_n$ is called a {\em descending path} and
its existence is denoted by $v_1\stackrel{*}{\to} v_n$.
We say that $v_n$ is a {\em descendant of $v_1$.}
We allow the empty path, $v \stackrel{*}{\to} v,$ for any $v$.
Consequently, $v\stackrel{*}{\to} w$ is a relation on $V$
which is reflexive and transitive but not necessarily symmetric.
We write $v \sim w$ if there is a finite path connecting
$v$ and $w$. (The edges of this path may have arbitrary
directions.)
The reduction rules $E$ are {\em (globally) confluent} if all
$v_1\sim v_2$ have a common descendant, i.e. $w\in V$ such that 
$v_1 \stackrel{*}{\to} w$ and $v_2 \stackrel{*}{\to} w$.
Finally, rules $E$ are {\em locally confluent} if for any
$v, w_1,w_2$ such that $v\to w_1,$ $v\to w_2,$ 
the elements $w_1,w_2$ have a common descendant.
Clearly, global confluence implies local confluence. 
However, the opposite implication fails, as shown in the following example.
The graph below contains infinitely many vertices and edges\vspace{.1in}:

\centerline{\parbox{2in}{\psfig{figure=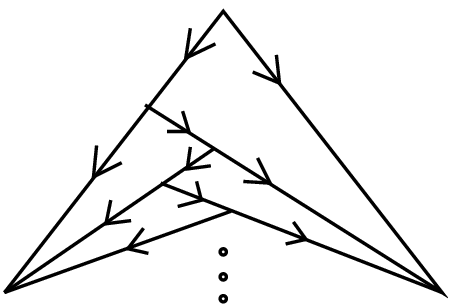,height=.9in}}
\vspace*{.1in}}

Nonetheless, under certain mild conditions on reduction rules,
local confluence implies global confluence.
We say that reduction rules are {\em terminal} if all descending paths are 
finite.\vspace*{.2in}

\noindent{\bf Diamond Lemma}\, \cite[Thm. 3]{Ne}
{\em If reduction rules are terminal then local confluence implies
global confluence.}\vspace*{.2in}

An example of an application of Diamond Lemma is the Jordan-H\"older
 theorem, which follows directly from this result.
Other applications of Diamond Lemma to ring theory and group theory 
are discussed in \cite{Be} and \cite{Sim}. Furthermore, Diamond Lemma
and the notion of confluence is used in mathematical logic:
in Church calculus, \cite{La, Ne}, in lambda calculus,
 \cite{BaN,Cu,Oh,La, Mi}, and in equational logic, \cite{OD}.
Additionally, it appears in computer science, in the theory of rewriting
systems and in the study of graph grammars, \cite{Eh1,Eh2,Na}.

%
\subsection{Confluence of linear objects}
\label{ss_lin_conf}
%

For our applications we need a  generalization of the notion of confluence to 
linear objects. For a ring $R,$ an $R$-linear reduction rule on a set 
$V$ is a pair $S: v\to \sum_{i=1}^n r_iv_i$ where $v,v_1,...,v_n\in V$ and
$r_1,...,r_n\in R$. Denote the free $R$-module over $V$ by $RV$.
For $X,Y\in RV,$ we write $X \stackrel{S}{\to} Y$ if $v$ appears with a
non-zero coefficient in $X$ and $Y$ is obtained from 
$X$ by replacing $v$ by $\sum_{i=1}^n r_iv_i$.
Finally, given a family of reduction rules, $\{S_i\}_{i\in I},$
we write $X \stackrel{*}{\to} Y$ if there is a sequence of
reduction rules leading from $X$ to $Y$.
Denote the $R$-submodule of $RV$ generated by $X-Y$ for all 
$X \stackrel{S_i}{\to} Y$ by $\R(S_i, i\in I)$ and write 
$X_1\sim X_2$ if $X_1-X_2\in \R(S_i, i\in I)$.
As before, we say that rules $\{S_i\}_{i\in I}$ are
{\em (globally) confluent} if any $X_1\sim X_2$ have a common descendant,
i.e. there is $Y$ such that $X_1 \stackrel{*}{\to} Y$ and
$X_2 \stackrel{*}{\to} Y$.
Finally, rules $\{S_i\}_{i\in I}$ are {\em locally confluent} on $V$
(respectively: on $RV$) if for any $X\in V$ (respectively: any $X\in RV$)
and any $Y_1,Y_2\in RV$ such that $X \stackrel{S_i}{\to} Y_1,$ 
$X \stackrel{S_j}{\to} Y_2,$ $Y_1$ and $Y_2$ have a common descendant.
Clearly, global confluence implies local confluence on $RV$, and terminal
local confluence on $RV$ implies global confluence.
However local confluence on $V$ does not imply local 
confluence on $RV$!
For example, let $V=\{v_1,v_2\}$ and let $S_1:v_1\to v_1+2v_2,$ 
$S_2: v_2\to v_2+2v_1$. Obviously, $S_1,S_2$ are locally confluent
on $V,$ since for no $v\in V$, $S_1(v)$ and $S_2(v)$ are simultaneously defined.
However $S_1$ and $S_2$ are not confluent on ${\mathbb R}V$! 

\begin{lemma}\label{no_descendant} 
$S_1(v_1+\sqrt{2}v_2)$ and $S_2(v_1+\sqrt{2}v_2)$
have no common descendant in ${\mathbb R}V$.
\end{lemma}

\begin{proof}
Notice that $S_i$ sends $a_1v_1+a_2v_2$ to $b_1v_1+b_2v_2,$ 
where $$\left({b_1 \atop b_2}\right)=M_i \left({a_1 \atop a_2}\right)
\quad \text{and}\quad M_1= \left(
\begin{array}{cc} 1 & 0\\ 2 & 1\\ \end{array}\right),
\quad M_2= \left(
\begin{array}{cc} 1 & 2\\ 0 & 1\\ \end{array}\right).$$
If $S_1(v_1+\sqrt{2}v_2)$ and $S_2(v_1+\sqrt{2}v_2)$ have a common descendant
$c_1v_1+c_2v_2$ for some $c_1,c_2\in \mathbb R$, 
then for certain products $N_1,N_2$ 
of matrices $M_1,M_2,$
$$N_1M_1\left({1 \atop \sqrt{2}}\right)=N_2M_2\left({1 \atop \sqrt{2}}\right)=
\left({c_1 \atop c_2}\right).$$
Irrationality of $\sqrt{2}$ implies $N_1M_1=N_2M_2$ as matrices in $SL(2,\Z)$.
However, $M_1,M_2$ generate a free semigroup in $SL(2,\Z)$.
Therefore, $N_1M_1\ne N_2M_2$ for any $N_1,N_2$.
\end{proof}

Nonetheless, we have

\begin{theorem}[Linear Diamond Lemma] 
Let $\{V_j\}_{j\in J}$ be a family of subsets of $V,$
such that $J$ is a well ordered set and 
$V_j\subset V_{j'}$ for $j<j'$ and $\bigcup_{j\in J} V_j=V$. 
Let $deg(v)=min\, \{j: v\in V_j\}$.
Consider a family of linear reduction rules on $V$ 
such that each of them sends an element of $V$ to a linear combination of
elements of smaller degree. Then\\
(1) these reduction rules are terminal,\\
(2) if this family is locally confluent on $V$ then 
it is also locally confluent on $RV$.
Therefore, by Diamond Lemma, such family of reduction rules is globally 
confluent on $RV$.
\end{theorem}

\begin{proof}
(1) Assume that there is an infinite chain $X_1\to X_2\to X_3\to ....$
Let $X_i=\sum_j^{n_i} c_{ij}v_{ij}$ and let $d_{ik}$ denote the 
$k$-th highest degree among degrees of $v_{i,1},...,v_{i,n_i}$.
Since $d_{11}\geq d_{21}\geq d_{31}\geq ...,$ the sequence stabilizes at 
certain place, which we denote by $N_1$. In other words $d_{k,1}=d_{N_1,1},$
for all $k\geq N_1$. Let $e_1=d_{N_1,1}$.
Similarly, $d_{N_1,2}\geq d_{N_1+1,2}\geq d_{N_1+2,2}\geq ...,$ 
stabilizes, let us say, at $N_2$-th place. Let $e_2=d_{N_2,2}$.
By continuing this process, we construct $e_1\geq e_2\geq e_3\geq ...$
This sequence stabilizes at some point as well -- let us say at $s$.
Then for any $k\geq N_s,$ the elements of $V$ appearing in 
$X_k=\sum_j^{n_i} c_{kj}v_{kj}$ have degrees $e_1,...,e_s$ 
(each of them may be appearing many times).
This, however, leads to contradiction since any reduction transformation 
replaces some $v$ by a linear combination of elements of 
$V$ of lower degree.\\
(2) Let $S_1: v_1\to \sum_{i=1}^{n_1} b_iw_i$ and $S_2: v_2\to 
\sum_{i=1}^{n_2} c_iz_i$.
Assume that $deg(v_1)<deg(v_2).$
We need to prove that for any $X,$
$S_1(X), S_2(X)$ have a common descendant.
Let $X=a_1v_1+a_2v_2+ X',$ where $X'$ is a linear combination of elements
of $V\setminus \{v_1,v_2\}.$
Since degrees of $w_1,..., w_{n_1}$ are smaller than that of $v_2,$
the elements $w_1,...,w_{n_1}$ are different than $v_2.$ If, additionally, 
$z_1,...,z_{n_2}\ne v_1$ then 
\begin{equation}\label{des}
S_2S_1(X)=a_1\sum_{i=1}^{n_1} b_iw_i+a_2\sum_{i=1}^{n_2} c_iz_i +X'=
S_1S_2(X)
\end{equation}
is a common descendant of $S_1(X)$ and $S_2(X)$ and the proof is complete. 
Therefore, assume now that one of the $z_i$'s, say $z_1$ is equal to $v_1$. 
If $a_2c_1=0$ then $S_1S_2(X)=S_2S_1(X)$ again. However, this may not be
the case if $a_2c_1\ne 0,$ since then
$$S_1S_2(X)=(a_1+a_2c_1)\sum_{i=1}^{n_1} b_iw_i+a_2\sum_{i=2}^{n_2} c_iz_i +
X'$$ and $S_2S_1(X)$ is as in (\ref{des}). Now, however,
$S_1S_2(X)=S_1S_2S_1(X)$ is a common descendant of
$S_1(X)$ and $S_2(X)$. 
\end{proof}

$X\in V$ is {\em irreducible} with respect to a given set of
reduction rules if none of these rules applies to $X$.
Denote the set of irreducible elements by $V_{irr}$.
Note that if $\{S_i\}_{i\in I}$ are terminal then 
$RV/{\mathcal R}(S_i, i\in I)$ is spanned by $V_{irr}$.
The opposite implication does not hold in general.

The combination of confluence and termination is a very strong
property of reduction rules.

\begin{theorem}\label{main}
(1) For any terminal rules $\{S_i\}_{i\in I}$ for $RV$
the following conditions are equivalent:
\begin{enumerate}
\item[(a)] $S_i,$ $i\in I,$ are locally confluent in $RV;$
\item[(b)] $S_i,$ $i\in I,$ are confluent in $RV;$
\item[(c)] For any $x\in RV$ there is a unique element 
$\psi(x) \in RV_{irr}$ such that $x\stackrel{*}{\to} \psi(x)$. 
\end{enumerate}
(2) If any of the above conditions holds then
$\psi: RV\to RV_{irr}$ is an $R$-linear map which factors
to an isomorphism $$\bar \psi: RV/{\mathcal R}(S_i, i\in I)\to
RV_{irr}.$$ 

Furthermore, $\psi$ is the identity
on $RV_{irr}$ and, consequently, $V_{irr}$ is a basis of 
$RV/{\mathcal R}(S_i, i\in I)$.
\end{theorem}

\begin{proof} 
(a) $\Rightarrow$ (b) by the Diamond Lemma.

(b) $\Rightarrow$ (c): Since the reduction rules are terminal, every $x\in RV$ 
has a descendant $y\in RV_{irr}$.
By confluence, $y$ is unique -- indeed, if 
$x \stackrel{*}{\to} y'\ne y$ and $y'\in RV_{irr}$ then $y\sim y'$ but 
they have no common descendants, contradicting the confluence assumption.

(c) $\Rightarrow$ (a) is obvious. 

(c) $\Rightarrow$ (2):
If $x\stackrel{*}{\to} y$ then $\psi(x)=\psi(y)$.
Since the relation $\sim$ defined at the beginning of Section 
\ref{ss_lin_conf} is the smallest equivalence relation
on $RV$ generated by $\stackrel{*}{\to},$ $x\sim y$ implies that
$\psi(x)=\psi(y)$. 
Therefore $\psi$ factors to 
$$\bar \psi: RV/{\mathcal R}(S_i, i\in I)=
RV/\sim\ \to RV_{irr}.$$

If we denote the ``obvious'' map 
$RV_{irr}\to RV\to RV/{\mathcal R}(S_i, i\in I)$ by $\imath$ then
clearly both $\imath \psi$ and $\psi \imath$ are identities on their
respective domains. Therefore $\bar \psi$ is a bijection and an 
$R$-linear map. Finally, $\psi$ is also $R$-linear, since it is a 
composition of linear maps $$RV\to RV/{\mathcal R}(S_i, i\in I)
\stackrel{\bar \psi}{\to} RV_{irr}.$$
\end{proof}

Only a few interesting terminal and confluent reduction 
systems on sets of graphs are known. Most of them 
appear in the context of representation theory of Lie
algebras of rank $\leq 2$ and
of associated quantum groups, c.f. Sections \ref{s_a1}-\ref{s_g2}.
See Section \ref{s_dichrom} for other examples.

%
\section{Graphs}\label{s_graphs}
%

%
\subsection{Abstract Graphs}\label{ss_abstract_graphs}
%

In a most general setting, a {\em labeled graph} is 
$\Gamma=(\V,\E,t,\tau,\Lambda,
\lambda,\nu),$ where $\V$ is a vertex set,
$\E$ is the set of {\em edge directions}, $t: \E\to \V$ is the {\em tail} map, 
$\tau: \E\to \E$ is the {\em change of direction
involution} which is fixed-point free. $\Lambda$ is a set of 
labels and $\lambda:\E\to \Lambda$ is a labeling function.
$\nu: \Lambda\to \Lambda$ is an involution such that $\nu \lambda=
\lambda\tau :\E\to \Lambda$. The function $t\tau: E \to V$ is 
called the {\em head map.}

An edge is a two-element set $\{e,\tau(e)\}$. 
The {\em valency} of $v\in V$ is the number of edge directions $e$ such that 
$t(e)=v$.
As mentioned in Introduction, we will sometimes specify a set 
of $1$-valent vertices
$V_{ext}(\Gamma)\subset V(\Gamma),$ called {\em external vertices}, and 
consider it as part of graph structure of $\Gamma$. 
The remaining vertices, $V_{int}(\Gamma)=V(\Gamma)\setminus V_{ext}(\Gamma),$ 
are {\em internal}.

Most definitions of graphs can be deduced from this one.
For example, a {\em partially directed} graph is
$\Gamma=(\V,\E,t,\tau,\emptyset,\Lambda,\lambda,\nu),$
such that $\Lambda=\{\pm 1,0\}$ and $\nu(x)=-x$. An edge 
$\{e,\tau(e)\}$ with $\lambda(e)=0$ is {\em undirected.}
Otherwise, its {\em direction} is either $e$ or $\tau(e)$ depending on
whether $\lambda(e)=1$ or $-1$. 

The reason for using edge directions, instead of
edges, is that in representation theory one considers graphs
whose edges are labeled by representations and have no canonical orientation.
If an edge direction is labeled by a representation V
then the opposite direction is labeled by the dual of $V$.

An embedding of $\Gamma_1$ into $\Gamma_2$ is\\ 
(1) a map $f: V_1\to V_2$ which is an embedding of internal vertices of 
$\Gamma_1$ into internal vertices of $\Gamma_2$,\\
(2) a map $g: \E_1\hookrightarrow \E_2,$
such that $t_2g=ft_1,$ $\tau_2g=g\tau_1$, and $g$ restricted to
$\{e\in \E_1: t(e) \text{\ is an internal vertex}\}$ is an
embedding.\\
(3) an embedding $h: \Lambda_1\hookrightarrow \Lambda_2$
such that $\lambda_2g=h\lambda_1$ and $\nu_2 h=h\nu_1$.

For any embedding $f :\Gamma\hookrightarrow \Gamma'$
and $e\in E(\Gamma'),$ neighborhood of $f^{-1}(e)$
has one of the following forms\footnote{To be precise, one
considers the topological realization of $\Gamma$ and the topological
neighborhood of $f^{-1}(e)$.}:

$$\begin{array}{ccccccccc}
\parbox{.35in}{\psfig{figure=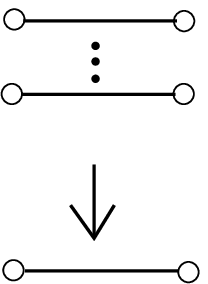,height=.6in}} & &
\parbox{.7in}{\psfig{figure=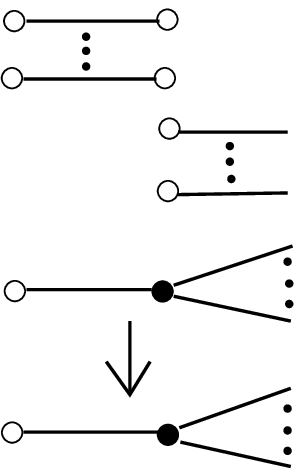,height=1in}} & &
\parbox{.7in}{\psfig{figure=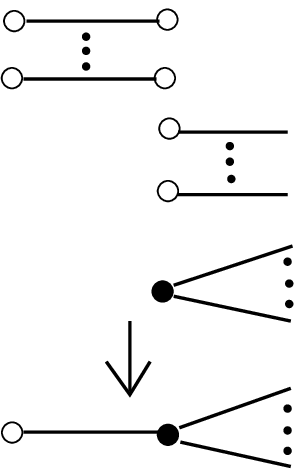,height=1in}} & &
\parbox{.7in}{\psfig{figure=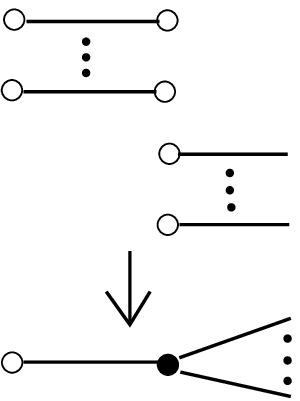,height=1in}} & &
\parbox{1in}{\psfig{figure=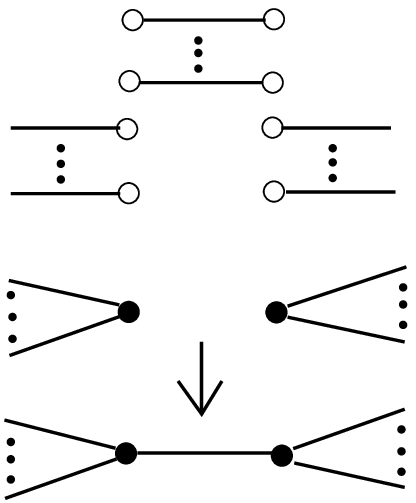,height=1in}}\\
 & & & & & & & & \\
\text{Type:\ } EE & & EI_1 & & EI_2 & &
EI_3 & & II_1\\
\end{array}$$

$$\begin{array}{ccccc}
\parbox{.8in}{\psfig{figure=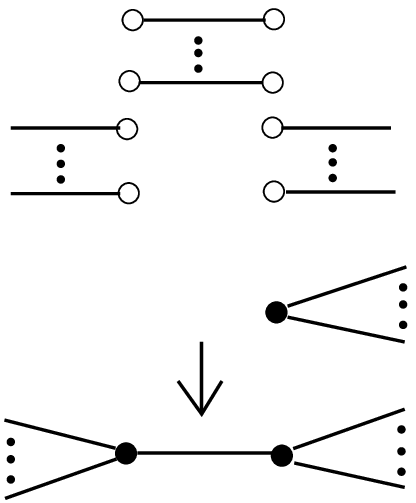,height=1in}} &
\parbox{.8in}{\psfig{figure=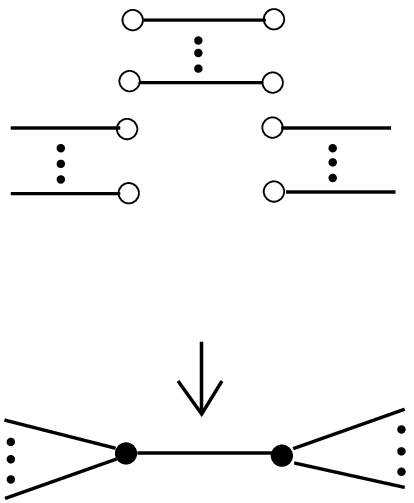,height=1in}} &
\parbox{.8in}{\psfig{figure=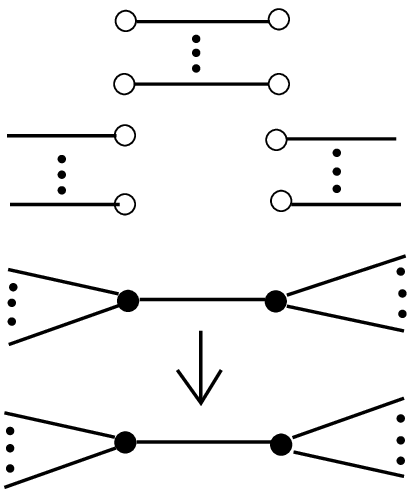,height=1in}} &
\parbox{.8in}{\psfig{figure=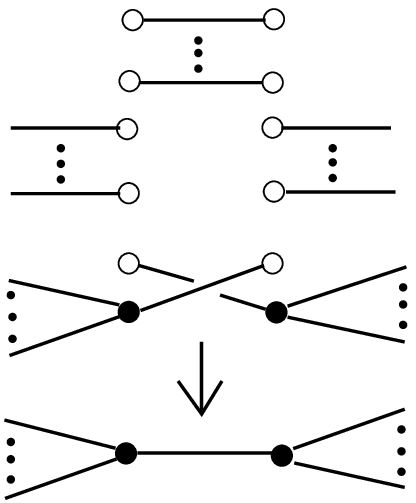,height=1in}} &
\parbox{.8in}{\psfig{figure=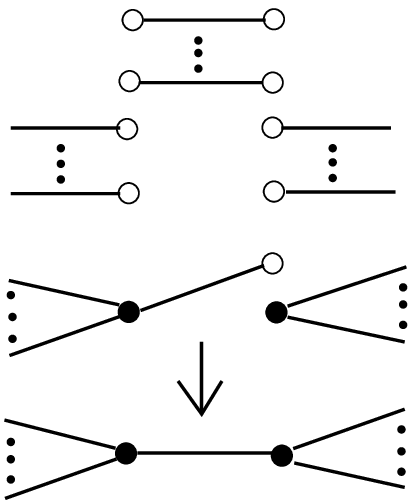,height=1in}}\\
 & & & &  \\
II_2  & II_3  & II_4  & II_5 & II_6\\
\end{array}$$

Above, black dots denote internal vertices and white dots
the external ones. Triple dots denote several parallel copies (possibly zero).
For the purpose of this classification we ignore edge directions.

%
\subsection{Graphs in manifolds}\label{ss_mfld_graphs}
%

Throughout the paper all manifolds are smooth.
A graph in a manifold $M,$ or {\em manifold graph}, 
is a subspace $\Gamma\subset M$ which looks locally 
like $1$-dimensional submanifold of $M$ (possibly with boundary) 
except for {\em internal vertices}, \diag{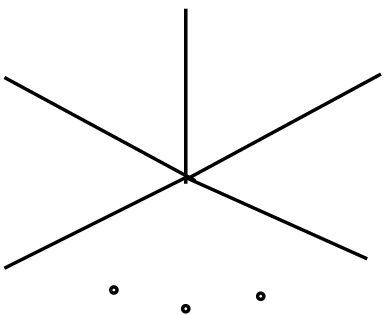}{.4in}\quad .
For simplicity, we do not allow $2$-valent vertices.

The points of $\Gamma\cap \p M=V_{ext}(\Gamma)$ are called 
{\em external vertices}.
Manifold graphs are considered up to isotopy of $M$ fixing $\p M$.
We denote the set of all vertices of $\Gamma$ by $V(\Gamma)$.
Note that $\Gamma\setminus V(\Gamma)$ is composed of open intervals
and circles (also called loops).

A manifold graph $\Gamma\subset M$ is {\em labeled} if there is specified
a set $\Lambda$ with an involution $\tau:\Lambda\to \Lambda$ and 
a labeling function 
$$\lambda:\left\{\text{orientations of connected components of 
$\Gamma\setminus V(\Gamma)$}\right\} \to \Lambda.$$ 
We require that if $o,\bar o$ are opposite orientations of the same
edge or circle in $\Gamma\setminus V(\Gamma)$ then $\lambda(\bar o)=
\tau(\lambda(o))$. 
A graph labeled by $\Lambda=\{0,\pm 1\},$ with $\tau(x)=-x,$ 
is called {\em partially oriented}. An edge or circle $e$ of $\Gamma$
is unoriented if $\lambda(e)=0$ and oriented otherwise.
Its orientation is the one labeled by $1$.

Note that if $M$ is connected, simply-connected, has connected boundary,
and $dim\, M\geq 4$ then
each graph in $M$ can be thought as a geometric realization of an abstract
graph. Such abstract graph is unique up to an insertion or deletion 
of $2$-valent internal vertices into edges or from edges.

An {\em embedding} of manifold graph $\Gamma_1\subset M_1$ into
$\Gamma_2\subset M_2$ is an embedding $f: M_1\hookrightarrow M_2$
of manifolds of equal dimensions,
which embeds a certain representative $\overline \Gamma_1$
of the isotopy class of $\Gamma_1\subset M_1$ into a certain representative 
$\overline \Gamma_2$ of the isotopy class of $\Gamma_2\subset M_2$ 
such that $f$ restricted to $\Gamma_1\setminus V_{ext}(\Gamma_1)$ 
is an open map into $\overline\Gamma_2$. We identify isotopic embeddings.
This definition implies that edges of $\Gamma_1$ are
mapped either into edges or circles of $\Gamma_2$. 

If $\Gamma_1\subset M_1,$ $\Gamma_2\subset M_2$ are labeled manifold graphs, 
then an embedding of $\Gamma_1$ into $\Gamma_2$ consists of a map 
$f: M_1\to M_2$ as above together with an embedding of the set of labels
$\imath: \Lambda_1\hookrightarrow \Lambda_2$ such that $f$ maps every edge 
or circle with some orientation, $e_1,$ of $\Gamma_1$ into
an edge or circle of $\Gamma_2,$ denoted by $e_2$ with coinciding orientation
such that $\lambda_2(e_2)=\imath \lambda_1(e_1)$.

\begin{example}
A graph embedding:
$\diag{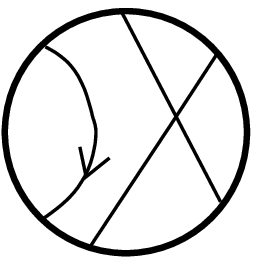}{.4in}\quad \hookrightarrow \quad 
\parbox{1in}{\psfig{figure=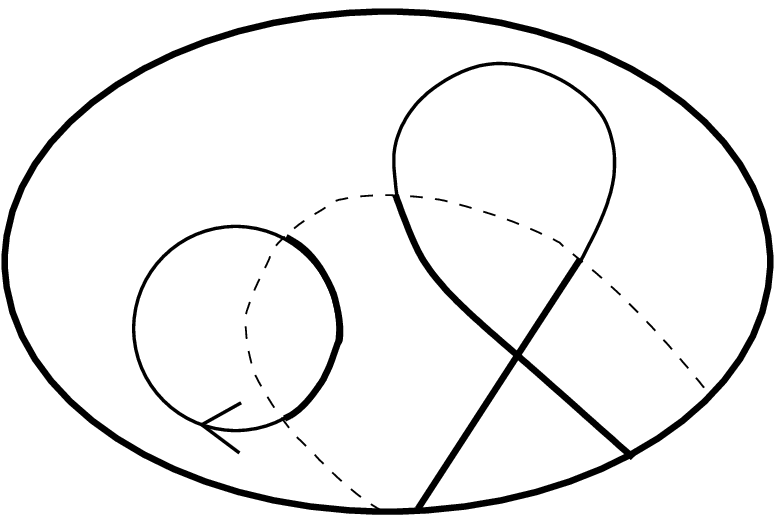,height=.6in}},$
and two non-embeddings:
$$\parbox{1in}{\psfig{figure=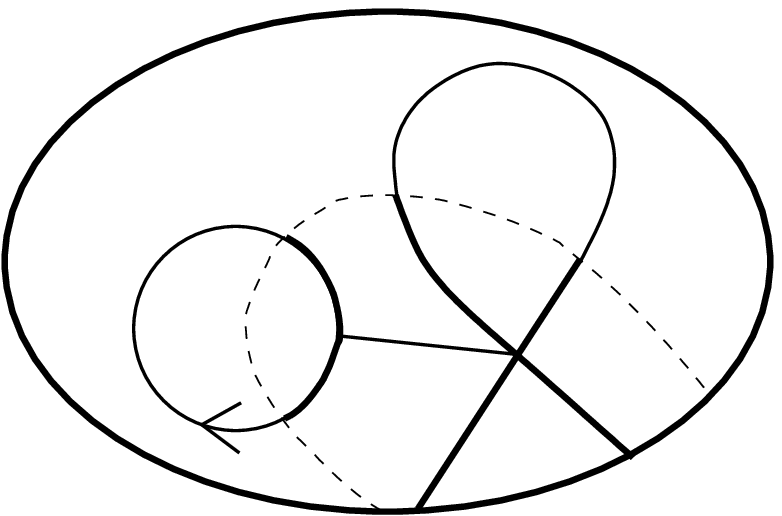,height=.6in}},\
\parbox{1in}{\psfig{figure=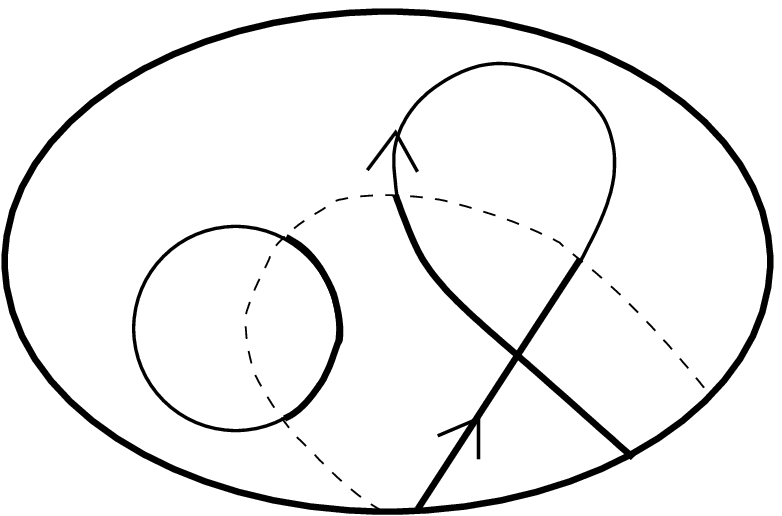,height=.6in}}$$
The above graph embedding is isotopic and, hence, identified with
the embedding \parbox{1in}{\psfig{figure=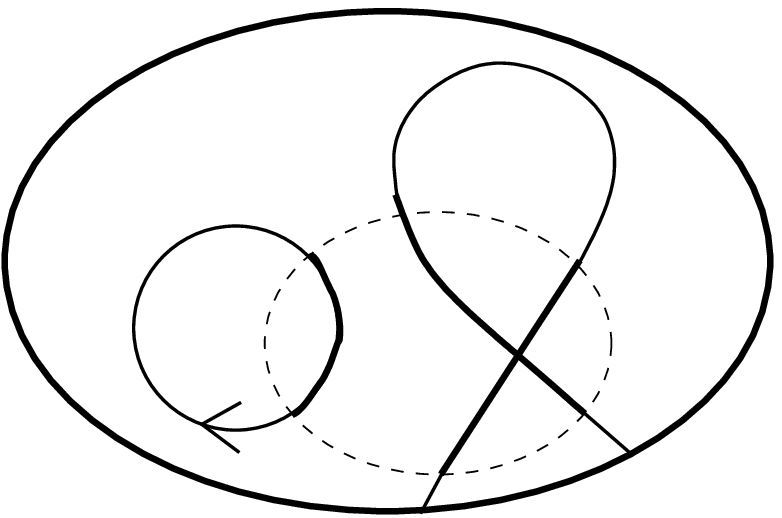,height=.6in}}.
\end{example}

The theory presented in this paper comes in two flavors: 
oriented and unoriented. In the first case all manifolds are 
oriented and all embeddings preserve orientations of manifolds. In the 
latter case, orientations of manifolds do not play any role. In both cases,
labelings of edge and circle orientations are preserved.
We will stress the difference between oriented and unoriented case 
whenever necessary, for example in Section \ref{s_a2}.

%
\subsection{Linear graphs}
\label{ss_linear}
%

Let $R$ be a ring. An {\em $R$-linear graph} is a formal
$R$-linear combination of graphs 
$\Gamma=\sum_{i=1}^k r_i \Gamma_i,$ together with specified bijections 
$V_{ext}(\Gamma_1)\simeq V_{ext}(\Gamma_2)\simeq .... \simeq V_{ext}(\Gamma_k)$
such that the corresponding external edge directions have identical 
labels.
(Since each external vertex is $1$-valent, its adjacent {\em external edge
direction} is well defined.)

An {\em $R$-linear manifold graph} in $M$ 
is a formal linear combination $\Gamma=\sum_{i=1}^n r_i \Gamma_i$ of graphs in 
$M$ such that their external vertices coincide
and the outward orientations of the corresponding external edges
have identical labels. These graphs are considered up to isotopy of 
$M$ fixing $\p M$.
For example,
$$T= \diag{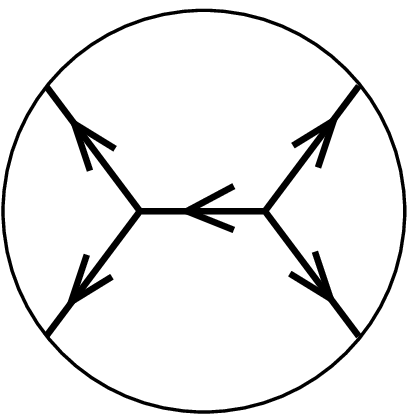}{.45in}-\diag{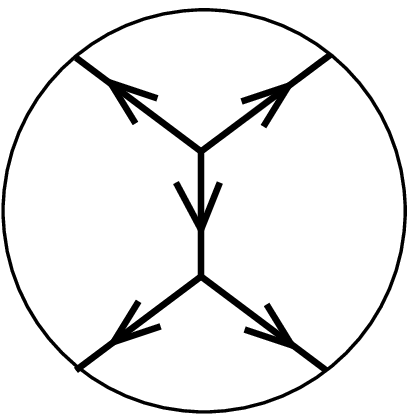}{.45in}$$
is a non-zero linear graph in $D^2$.

%
\subsection{Reduction rules on graphs}
\label{ss_reduction}
%

We are going to apply the theory of confluence to the problems stated in
Introduction, with a particular focus on graphs arising as Kuperberg's spiders
webs. Certain versions of this method were used implicitly already in 
\cite{Ja,Ku-spider,Ku-G2} and \cite{Ye}. Nonetheless, to our knowledge
the subtle difference between local confluence on $V$ and on $RV$ discussed 
in Sec. \ref{ss_lin_conf} has never been observed.

For our purposes, the set of objects, $V,$ considered in 
Section \ref{s_confluence} is either a set of abstract graphs or a set of
graphs in a given manifold $M$. 

In the first case, for a given ring $R,$ a {\em graph reduction rule} is a 
pair denoted by $T_0 \stackrel{S}{\to} T,$ where
$T_0$ is a graph, $T=\sum_{i=1}^k r_i T_i$ 
is an $R$-linear graph and the external vertices
of $T_0$ and $T_i$'s are identified via a bijection 
such that the corresponding external edge directions
have identical labels.
Any graph reduction rule $T_0 \stackrel{S}{\to} T$ defines reductions
of graphs $\Gamma$ as follows:
For any embedding $\Gamma_0\hookrightarrow \Gamma$
we obtain new graphs $\Gamma_i$ by replacing $T_0$ in $\Gamma$ by $T_i$.
We say that reduction $\Gamma_0\to \sum_{i=1}^k r_i \Gamma_i$ is induced
by $T_0 \stackrel{S}{\to} T$ and we denote that fact by putting $S$ above the
arrow, $\Gamma\stackrel{S}{\to} \sum_{i=1}^k r_i\Gamma_i$. 
Therefore, unlike in Section \ref{s_confluence}, we use one symbol 
(here, $S$) to denote many reduction rules
arising from $T_0 \stackrel{S}{\to} T$.

Similarly, a {\em graph reduction rule} for manifold graphs is a pair 
$T_0 \stackrel{S}{\to} T,$ where
$T_0$ is a graph in some manifold $M_0$ and $T=\sum_{i=1}^k r_i T_i$ 
is an $R$-linear graph in the same manifold such that the external vertices
of $T_0$ and $T_i$'s coincide and the corresponding external edge orientations
have identical labels.
For any embedding $M_0\hookrightarrow M$ and
a graph $\Gamma\subset M$ such that $\Gamma\cap M_0=T_0,$
we obtain new graphs $\Gamma_i\subset M$ by replacing $T_0$ in 
$\Gamma$ by $T_i$.
We say that reduction $\Gamma_0\to \sum_{i=1}^k r_i \Gamma_i$ is induced
by $T_0 \stackrel{S}{\to} T$ and we denote that fact by 
$\Gamma\stackrel{S}{\to} \sum_{i=1}^k r_i\Gamma_i$. 

By analogy to the notation in Section \ref{s_confluence}, 
we use ${\mathcal R}(S_i, i\in I)\subset RV$ to denote
the submodule generated by all linear graphs 
$\Gamma - \sum_{i=1}^k r_i\Gamma_i$ coming from graph reductions
$\Gamma \stackrel{S_i}{\to}\sum_{i=1}^k r_i\Gamma_i,$ for $i\in I.$

%
\subsection{Proving confluence of reduction rules of abstract graphs}
\label{ss_abs_overlaps}
%
An {\em overlap} of graphs $\Gamma_1$ and $\Gamma_2$ is a graph $\Gamma$ and
pair of graph embeddings $(\Gamma_1\hookrightarrow \Gamma,
\Gamma_2\hookrightarrow \Gamma)$.
If $\{T_{i0}\stackrel{S_i}{\to} \sum_k r_{ik}T_{ik}\}_{i\in I}$ 
is a set of reduction rules of abstract graphs with coefficients in $R,$
then each overlap
$O=(\imath_1:T_{i0}\hookrightarrow \Gamma, 
\imath_2:T_{j0} \hookrightarrow \Gamma)$ 
leads to two different reductions of $\Gamma$.
We say that reduction rules $\{S_i\}_{i\in I}$ are {\em locally confluent on}
$O=(\Gamma_1\hookrightarrow \Gamma, \Gamma_2\hookrightarrow \Gamma)$
if for any $i,j$ such that $\Gamma_{i0}=\Gamma_1,$
$\Gamma_{j0}=\Gamma_2,$ the two reductions 
$\Gamma\stackrel{S_i}{\to} \sum_k r_{ik}\Gamma_k$ 
and $\Gamma\stackrel{S_j}{\to} \sum_k r_{jk} \Gamma_k',$ arising from 
this overlap have a common descendant.

Note that $\{S_i\}_{i\in I}$ are locally confluent on $V$ if and only
if they are locally confluent on all overlaps of graphs $T_{i0},$ $i\in I$.

We say that $O=(\imath_1:\Gamma_1\hookrightarrow \Gamma, 
\imath_2:\Gamma_2 \hookrightarrow \Gamma)$ factors through
$O'=(\imath_1':\Gamma_1\hookrightarrow \Gamma', 
\imath_2':\Gamma_2 \hookrightarrow \Gamma')$ 
if there is an embedding 
$f: \Gamma'\hookrightarrow \Gamma$,
such that $\imath_1=f\imath_1',$ $\imath_2=f\imath_2'$.
If reduction rules $\{S_i\}_{i\in I}$ are locally confluent on
$O$ then they are locally confluent on all overlaps which factor through $O$.
An overlap with no factorizations other than the identity is {\em irreducible}.

We are going to show that the following types of factorizations 
$(\imath_1':\Gamma_1\hookrightarrow \Gamma', 
\imath_2':\Gamma_2 \hookrightarrow \Gamma')\stackrel{f}{\to}
(\imath_1:\Gamma_1\hookrightarrow \Gamma, 
\imath_2:\Gamma_2 \hookrightarrow \Gamma)$
reduce every overlap to an irreducible one:

\begin{enumerate}
\item If $V(\Gamma_1)\cup V(\Gamma_2)$ is a proper subset of $V(\Gamma)$
then let $\Gamma'$ be a graph obtained from $\Gamma$ by 
removing vertices in $V(\Gamma)\setminus (V(\Gamma_1)\cup V(\Gamma_2)).$
Let $\imath_1'=\imath_1,$ $\imath_2'=\imath_1,$ and let $f$ be the obvious
embedding.
\item If $V_{int}(\Gamma_1)\cup V_{int}(\Gamma_2)$ is a proper 
subset\footnote{Since graph embeddings send internal vertices to 
internal vertices, $V_{int}(\Gamma_1)\cup V_{int}(\Gamma_2)\subset 
V_{int}(\Gamma).$} of $V_{int}(\Gamma)$, then
let $\Gamma'$ be a graph obtained from $\Gamma$ by changing the internal 
vertices in $V_{int}(\Gamma)\setminus (V_{int}(\Gamma_1)\cup 
V_{int}(\Gamma_2))$ to external ones. Let $f$ be the obvious
embedding.
\item If $E(\Gamma_1)\cup E(\Gamma_2)\subset E(\Gamma)$ is a proper subset, 
then let $\Gamma'$ be $\Gamma$ with the edge directions in 
$E(\Gamma)\setminus (E(\Gamma_1)\cup E(\Gamma_2))$ removed. $f$ is the obvious
embedding.

\item Let $E_{ext}(\Gamma)$ denote 
$\{e\in E(\Gamma): t(e),h(e)\in V_{ext}(\Gamma)\}.$
If $e_1\in E_{ext}(\Gamma_i),$ $\imath_i(e_1)=\imath_j(e_2),$
$e_1\ne e_2,$ for some $i,j\in \{1,2\},$ then let $\Gamma'$ be $\Gamma$ 
with two extra 
external vertices $v_1,v_2$ and two new edge directions $e',\tau(e')$
forming an edge connecting $v_1$ and $v_2.$ Let $\imath_1',\imath_2'$ 
coincide with $\imath_1,\imath_2$, except for $\imath_i'$ sending 
$e_1,\tau(e_1)$ to $e',\tau(e')$ and sending $t(e_1),h(e_1)$ to $v_1,v_2.$
Let $f(v_1)=t(e),$ $f(v_2)=h(e),$ $f(e')=e,$ and let $f$ be the identity on
the remaining vertices and edges.

\item If $e_1\in E_{ext}(\Gamma_i)$ and there is no edge $e_2$ as in (4) but 
$\imath_i(t(e_1))=\imath_j(v)$, $t(e_1)\ne v,$
for some $i,j\in \{1,2\},$ then let $\Gamma'$ be 
$\Gamma\setminus \{\imath_i(e_1),\tau(\imath_i(e_1))\}$ with an extra 
external vertex $w$ and two new edge directions $e',\tau(e')$
forming an edge connecting $w$ and $\imath_i(h(e_1)).$ 
Let $\imath_1',\imath_2'$ 
coincide with $\imath_1,\imath_2$, except for $\imath_i'$ sending 
$e_1, \tau(e_1)$ to $e',\tau(e')$ and sending $t(e_1)$ to $w.$
Let $f(w)=\imath_i(t(e_1)),$ $f(e')=\imath_i(e_1),$ 
$f(\tau(e'))=\imath_i(\tau(e_1)),$ 
and let $f$ be the identity on the remaining vertices and edge directions.

\item If $\imath_1^{-1}(e)$ and $\imath_2^{-1}(e)$ are of type $II_5$
(as defined at the end of Section \ref{ss_abstract_graphs}) then
consider factorization\\

\noindent
\parbox{3in}{\psfig{figure=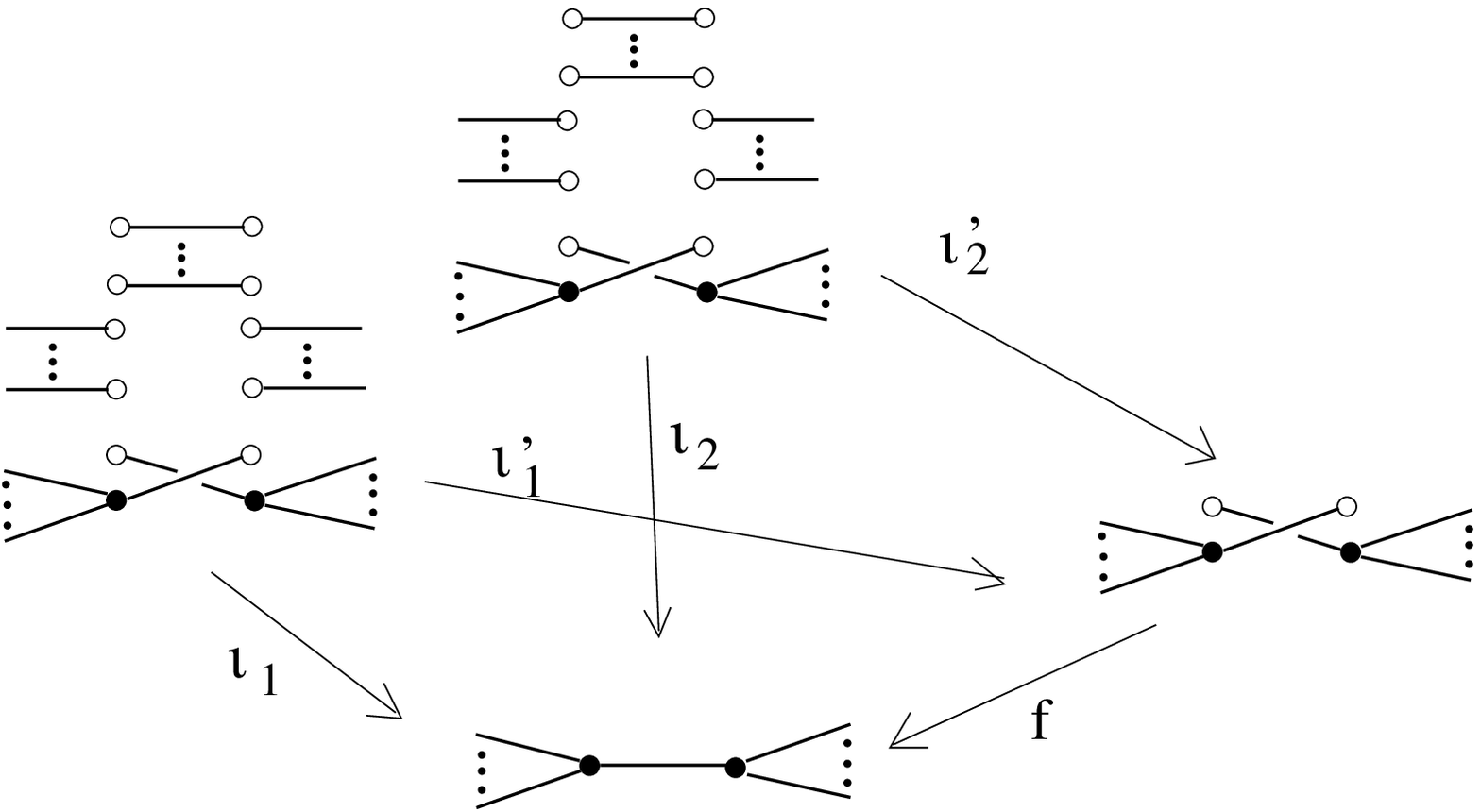,height=2in}}\vspace*{.3in}

\noindent (For simplicity the edge directions are ignored in the above 
picture.) Analogously, if $\imath_1^{-1}(e)$ and $\imath_2^{-1}(e)$ are 
(in some order) of 
types $(II_1,II_5),$ $(II_1,II_6),$ $(II_2,II_5),$ $(II_2,II_6),$ 
$(II_3,II_5),$ $(II_3,II_6),$ $(II_5,II_6),$ $(II_6,II_6)$
then we perform similar factorizations.
\end{enumerate}

\begin{theorem}\label{abstract_overlap}
An overlap $O=(\imath_1:\Gamma_1\hookrightarrow \Gamma, \imath_2:\Gamma_2 
\hookrightarrow \Gamma)$ is irreducible iff
\begin{enumerate}
\item $V(\Gamma)=V(\Gamma_1)\cup V(\Gamma_2).$
\item $V_{int}(\Gamma)=V_{int}(\Gamma_1)\cup V_{int}(\Gamma_2).$
\item $E(\Gamma)=E(\Gamma_1)\cup E(\Gamma_2)$
\item If $e_1\in E_{ext}(\Gamma_i)$ and $\imath_i(e_1)=\imath_j(e_2)$
then $j=i$ and $e_1=e_2.$
\item If $e\in E_{ext}(\Gamma_i),$ 
$\imath_i(t(e))=\imath_j(v)$, then $i=j$ and $t(e)= v.$
\item $\Gamma$ has no edges of types 
$(II_1,II_5),$ $(II_1,II_6),$ $(II_2,II_5),$ $(II_2,II_6),$ 
$(II_3,II_5),$ $(II_3,II_6),$ $(II_5,II_5),$ $(II_5,II_6),$ $(II_6,II_6).$
\end{enumerate}
\end{theorem}

\begin{proof} $\Rightarrow$
If one of the conditions does not hold then $O$ admits
one of the factorizations described above.

$\Leftarrow$ Suppose that $O$ satisfies (1)-(6) above and $$f:
(\imath_1':\Gamma_1\hookrightarrow \Gamma', 
\imath_2':\Gamma_2 \hookrightarrow \Gamma')\to
(\imath_1:\Gamma_1\hookrightarrow \Gamma, 
\imath_2:\Gamma_2 \hookrightarrow \Gamma)=O$$
is a factorization of $O.$

\begin{lemma} $f: V(\Gamma')\to V(\Gamma)$ is a bijection.
\end{lemma}

\begin{proof} $f(V(\Gamma'))\supset 
f\imath_1'V(\Gamma_1)\cup f\imath_2'V(\Gamma_2)=\imath_1V(\Gamma_1)\cup 
\imath_2V(\Gamma_2)$ is by (1) equal to $V(\Gamma).$ Therefore,
$f: V(\Gamma')\to V(\Gamma)$ is onto.

Suppose $f(v_1)=f(v_2)=v$ for $v_1\ne v_2,$ $v_1,v_2\in V(\Gamma').$
By definition of graph embedding, $f$ is $1$-$1$ on 
$V_{int}(\Gamma').$ Hence, at least one of the vertices $v_1,v_2$ is external, 
say $v_1.$ Consider two cases:

(a) $v$ is external. Then $v_2$ is external as well.
Let $e_1,e_2$ be edge directions with tails $v_1,v_2.$
Since $v$ is $1$-valent, $f(e_1)=f(e_2)\stackrel{def}{=} e\in E(\Gamma).$
Since $f(h(e_1))=f(h(e_2))=h(e),$ at least one of the vertices 
$h(e_1),h(e_2)$ is external, contradicting (4).

(b) $v$ is internal. By (2), we can choose $v_2$ to be an internal vertex.
Since $v_1$ has valency one, there is a unique vertex $w$ in $\Gamma'$
connected with $v$ by an edge $e'.$ By (5), $w$ is internal.
Let $e=f(e').$ Then the preimage, $f^{-1}$, of the neighborhood of $e$ 
must be of type $II_5$ or $II_6.$
Since $\imath_1,$ $\imath_2$ satisfy (1)-(4) and 
factorize through $f$, they must be of one of the types listed in (6).
Contradiction.
\end{proof}

\begin{corollary}\label{f_bij} $f:V_{int}(\Gamma')\to V_{int}(\Gamma)$
and $f:V_{ext}(\Gamma')\to V_{ext}(\Gamma)$
are bijections.
\end{corollary}

\begin{proof}
Since $f$ is $1$-$1$ it is enough to prove that
(a) $f(V_{int}(\Gamma'))= V_{int}(\Gamma)$ and (b)
$f(V_{ext}(\Gamma'))= V_{ext}(\Gamma)$.

\noindent (a) Since $f$ is a graph embedding, $f(V_{int}(\Gamma'))\subset 
V_{int}(\Gamma)$.
By (2), $\imath_1(V_{int}(\Gamma_1))\cup \imath_2(V_{int}(\Gamma_2))=
V_{int}(\Gamma)$, and since $\imath_1,\imath_2$ factor through $f,$
$f(V_{int}(\Gamma'))= V_{int}(\Gamma)$. 

\noindent (b) By previous lemma and by (a), 
$|V(\Gamma')|= |V(\Gamma)|,$ $|V_{int}(\Gamma')|= |V_{int}(\Gamma)|.$
Hence, $|V_{ext}(\Gamma')|= |V_{ext}(\Gamma)|$ and the statement
follows from the fact that $f$ is $1$-$1$.
\end{proof}

\begin{proposition}
$f: E(\Gamma')\to E(\Gamma)$ is a bijection.
\end{proposition}

\begin{proof}
By (3), $f$ is onto.
Suppose that $f(e_1)=f(e_2)=e,$ $e_1\ne e_2.$ Since $f$ is a bijection on
vertices, $t(e_1)=t(e_2),$ $h(e_1)=h(e_2).$ Since $t(e_1),$ $h(e_2)$ are
at least $2$ valent, they are internal and, consequently, 
$f$ maps two internal edge directions to a single edge direction 
and, therefore, it is not a graph embedding.
\end{proof} 

Therefore $f$ is the identity and the proof of Theorem \ref{abstract_overlap}
is completed.
\end{proof}

Since each overlap of $\Gamma_1$ and $\Gamma_2$ satisfying 
Theorem \ref{abstract_overlap}(1)-(3) is obtained as a quotient of
the disjoint union of $\Gamma_1$ and $\Gamma_2$, the number of such
overlaps is finite.

\begin{corollary}\label{finite_overlap}
Any two abstract graphs have a finite number of irreducible overlaps only.
\end{corollary}

Consider an overlap of $\Gamma_1$ and $\Gamma_2.$ By applying factorizations
of types (1)-(5), we obtain an overlap satisfying conditions (1)-(5) of
Theorem \ref{abstract_overlap}. Observe, that if an overlap $O$ satisfies
these conditions then for every factorization $f:O'\to O$ of type (6),
$O'$ satisfies (1)-(5) as well. Furthermore, observe that for each 
factorization $$f:
(\imath_1':\Gamma_1\hookrightarrow \Gamma', 
\imath_2':\Gamma_2 \hookrightarrow \Gamma')\to
(\imath_1:\Gamma_1\hookrightarrow \Gamma, 
\imath_2:\Gamma_2 \hookrightarrow \Gamma)=O$$
of type (6), either (a) the number of connected components of $\Gamma'$ is
larger than the number of components of $\Gamma$, or
(b) the number of cycles of $\Gamma'$ is
lower than the number of cycles of $\Gamma$.
Since the number of connected components is bounded above by the
(unchanging) number of vertices, every sequence of factorizations of type
(6) is finite.
Therefore, we proved:

\begin{corollary}\label{cor_overlap}
(1) Each overlap can be reduced to an irreducible one by a finite number
of factorizations of types (1)-(6).\\
(2) $\{T_{i0}\stackrel{S_i}{\to} T_i\}_{i\in I}$ are locally confluent 
if they are locally confluent on all irreducible overlaps of pairs of 
graphs in $\{T_{i0}\}_{i\in I}$.
\end{corollary}

Therefore, one has an effective procedure for deciding whether 
any finite set of reduction rules on graphs is locally 
confluent.

\begin{example}\label{overlap_example}
Graphs $\diag{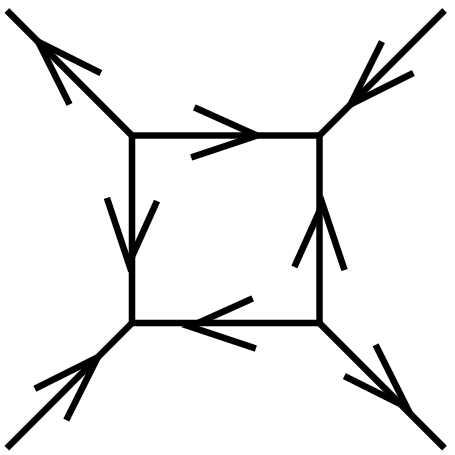}{.5in}$ and $\diag{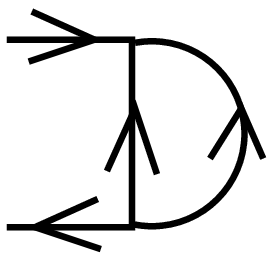}{.3in}$
have $5$ different irreducible overlaps:
four overlaps of the form
$$\diag{square}{.5in} \longrightarrow \diag{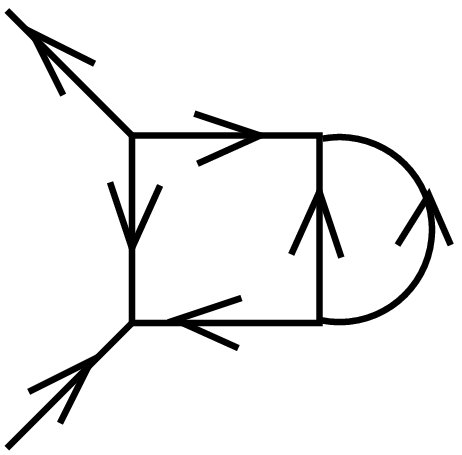}{.5in}
\longleftarrow \diag{olap2}{.3in}$$ 
and
one ``trivial'' overlap
$$\diag{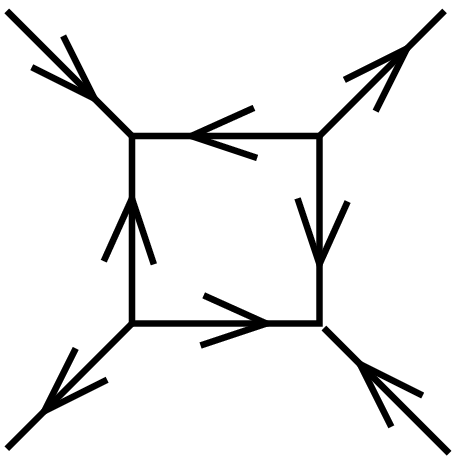}{.5in} \longrightarrow \diag{squareo}{.5in}\hspace*{.2in}
\diag{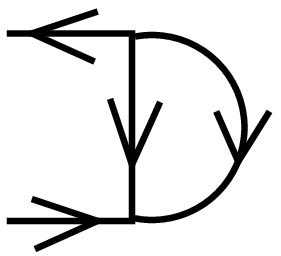}{.3in}\longleftarrow \diag{olap3}{.3in}$$
\end{example}

\begin{example} $\Gamma_1=\Gamma_2=$\diag{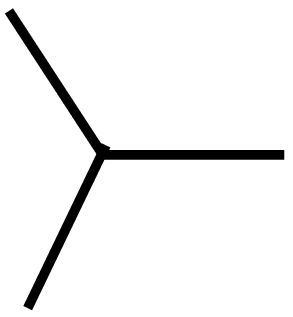}{.3in}
have three different irreducible overlaps of the form 
$$\diag{graphy}{.3in}\longrightarrow \diag{graphy}{.3in}\longleftarrow
\diag{graphy}{.3in}$$ and one trivial overlap,
$$\diag{graphy}{.3in}\longrightarrow \diag{graphy}{.3in}\hspace*{.2in}
\diag{graphy}{.3in} \longleftarrow \diag{graphy}{.3in}\hspace*{.2in}.$$
\end{example}

Notice that the embeddings of $\Gamma_1,\Gamma_2$ into 
\diag{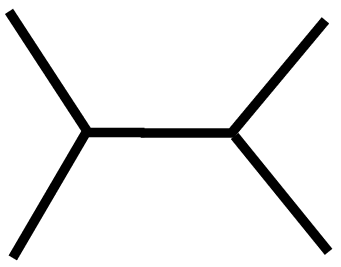}{.2in}\hspace*{.2in}, 
\diag{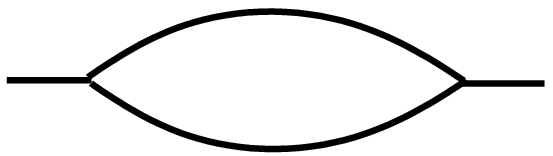}{.2in}\hspace*{.6in}, \diag{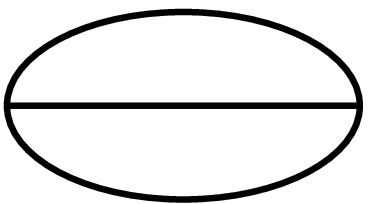}{.25in}\hspace*{.3in} 
are not irreducible since they factor through the trivial overlap.

%
\subsection{Proving confluence of reduction rules of surface graphs}
\label{ss_overlaps}
%

An {\em overlap} of manifold graphs $\Gamma_1\subset M_1,$ 
$\Gamma_2\subset M_2$ is a graph $\Gamma$ in a manifold $M$
together with isotopy classes of embeddings
$(\Gamma_1,M_1)\hookrightarrow (\Gamma,M),
(\Gamma_2,M_2)\hookrightarrow (\Gamma,M)$.\\
$f: O'=((\Gamma_1,M_1)\hookrightarrow (\Gamma',M'),
(\Gamma_2,M_2)\hookrightarrow (\Gamma',M'))\to$\\
\hspace*{1.5in} $((\Gamma_1,M_1)\hookrightarrow (\Gamma,M),
(\Gamma_2,M_2)\hookrightarrow (\Gamma,M))=O$\\
is a factorization of $O$, if for certain representatives 
$\imath_1':(\Gamma_1,M_1)\hookrightarrow (\Gamma,M),$
$\imath_2':(\Gamma_2,M_2)\hookrightarrow (\Gamma,M),$
of embeddings of $O',$
$f\imath_1,$ $f\imath_2$ belong to isotopy classes of embeddings of $O.$

As before, given reduction rules, 
$\{T_{i0}\stackrel{S_i}{\to} \sum_k r_{ik}T_{ik}\}_{i\in I},$
where $S_i$ takes place in a manifold $M_i,$ each overlap 
$O=(\imath_1:(T_{i0},M_i)\hookrightarrow (\Gamma,M),
\imath_2:(T_{j0},M_j)\hookrightarrow (\Gamma,M))$ 
leads to two different reductions of $\Gamma$.
Rules $\{S_i\}_{i\in I}$ are locally confluent if they are 
locally confluent on all overlaps of graphs $T_{i0},$ $i\in I$. 
As before we consider factorization of overlaps and observe that 
if rules $\{S_i\}_{i\in I}$ are locally confluent on $O$ then they are 
locally confluent on all overlaps which factor through $O$.

A factorization $f:O'\to O$ is trivial if 
$f: M'\to M$ is isotopic to a homeomorphism. As before, an overlap is {\em 
irreducible} if it does not admit a non-trivial factorization.

In this section we are going to develop an algorithm for proving local 
confluence of overlaps of surface graphs.
Observe that we cannot apply verbatim the method of the previous section
to our current setting since
Corollary \ref{finite_overlap} and Corollary 
\ref{cor_overlap}(1) and (2) fail for surface graphs:

\begin{lemma}\label{overlap_trouble}
(1) If every component of $F$ has a non-empty boundary then
no overlap $((\Gamma_1, F_1)\hookrightarrow 
(\Gamma, F), (\Gamma_2,F_2) \hookrightarrow (\Gamma,F))$
is irreducible.\\
(2) If every component of $F$ has a non-empty boundary then
no overlap $((\Gamma_1, F_1)\hookrightarrow 
(\Gamma, F), (\Gamma_2,F_2) \hookrightarrow (\Gamma,F))$
factors through an irreducible one.\\
(3) If $F$ is closed, $F\ne S^2, RP^2,$ and 
$\Gamma$ is either empty or it is a contractible loop in $F,$ then 
$(\Gamma,F)$ has infinitely many irreducible overlaps with itself.
\end{lemma}

\begin{proof}
(1) Let $F'$ be $F$ with a disk removed from one of its components, $C.$
By imagining the disk lying ``very close'' to $\partial C,$
one can isotope $\imath_1,\imath_2$ to $\imath_1', \imath_2'$ so that 
$\imath_1'(F_1)\cup \imath_2'(F_2)\subset F'$. Consequently,
$(\imath_1,\imath_2)$ factors through 
$(\imath_1': (\Gamma_1, F_1)\hookrightarrow 
(\Gamma',F'), \imath_2': (\Gamma_2,F_2) \hookrightarrow 
(\Gamma',F'))$ via the embedding $f: F'\to F$.
This is a non-trivial factorization, contradicting the initial assumption.

(2) If an overlap as above factors through
$(\imath_1: (\Gamma_1, F_1)\hookrightarrow (\tilde\Gamma, \tilde F), 
\imath_2: (\Gamma_2,F_2) \hookrightarrow 
(\tilde \Gamma,\tilde F))$ then $\tilde F\subset F$ and consequently, every
component of $\tilde F$ has a non-empty boundary.

(3) For any diffeomorphism $f:(\Gamma,F)\to (\Gamma,F),$ diffeomorphisms
$\imath_1=f: (\Gamma,F)\to (\Gamma,F)$ and the identity map 
$\imath_2=f: (\Gamma,F)\to (\Gamma,F)$ form a an irreducible overlap 
which we denote by $O_f$. Notice that $O_f = O_{f'}$ if and only 
if $f'$ is isotopic to $f$. Since the mapping class group of 
$F$ is infinite, there are infinitely
many irreducible overlaps of this type.
\end{proof}

We will attempt to resolve these difficulties now.
$\{O_j\}_{j\in J}$ is a {\em basis} of overlaps of
$(\Gamma_1,F_1)$ and $(\Gamma_2,F_2)$ if every overlap of these surface graphs
factors through $O_j$ for some $j\in J$.

\begin{corollary}
The rules $S_i:\Gamma_{i0}\to \sum_k r_{ik} \Gamma_{ik},$ $i\in I,$
are locally confluent, if they are locally confluent on a certain basis of
overlaps of pairs of graphs in $\{(\Gamma_{i0}, F_i)\}_{i\in I}$.
\end{corollary}

Lemma \ref{overlap_trouble}(2) shows that not every pair 
of graphs in surfaces has a finite basis of overlaps. Furthermore, basis
of overlaps are generally not unique.
Nonetheless, we are going to show that any two simple graphs in surfaces
have a finite basis of overlaps. We say that $\Gamma\subset F$ is 
{\em simple} if $\Gamma$ is connected and
every component $C$ of $F\setminus \Gamma$ is either $D^2$ or an annulus
whose one boundary component lies in $\Gamma$ and the other in $\p F$. 

\begin{theorem}\label{finite_basis}
Any two simple graphs have a finite basis of overlaps.
\end{theorem}

Our proof is also an algorithm for finding such a finite basis.

For any graph $\Gamma \hookrightarrow F$ there is an $\ve_0>0$ such that
$\ve$-neighborhoods of $\Gamma$ in $F$ are diffeomorphic
to each other for all $\ve<\ve_0$. Denote such $\ve$-neighborhood 
by $\nu(\Gamma)$ and we call it a framing of $\Gamma$. 
Each framing of $\Gamma$ retracts onto $\Gamma$ and
each finite abstract graph has finitely many different 
framings only. 

Any overlap of $(\Gamma_1,\nu(\Gamma_1))$ and
$(\Gamma_2,\nu(\Gamma_2))$ factors through an overlap
$\imath_1:(\Gamma_1,\nu(\Gamma_1))\hookrightarrow
(\Gamma,F),\imath_2:(\Gamma_2,\nu(\Gamma_2))\hookrightarrow
(\Gamma,F),$ such that $(\imath_1:\Gamma_1\hookrightarrow
\Gamma,\imath_2:\Gamma_2\hookrightarrow \Gamma)$ is an irreducible overlap 
of abstract graphs and $F$ is a framing of $\Gamma$.
Consequently, such overlaps form a basis of overlaps of 
$(\Gamma_1,\nu(\Gamma_1))$ and $(\Gamma_2,\nu(\Gamma_2))$.
Denote them by $O_1,...,O_d$.

Now assume that $\Gamma_1$ and $\Gamma_2$ are embedded into
$F_1,F_2$ in such way that they are simple graphs.
We extend each basic overlap $O_i=(\imath_1:(\Gamma_1,\nu(\Gamma_1))
\hookrightarrow (\Gamma,F),\imath_2:(\Gamma_2,\nu(\Gamma_2))\hookrightarrow
(\Gamma,F))$ constructed above to an overlap of 
$(\Gamma_1,F_1)$ and $(\Gamma_2,F_2)$ as follows:
Every component $B$ of $\p F,$ disjoint from $\Gamma,$ 
is parallel to a unique cycle $\alpha_B\subset \Gamma$.
If the preimage $\imath_i^{-1}(\alpha_B)$ for either $i=1$ or $2$ 
is a circle in $F_i$ which bounds a disk $D_i\subset F_i$ containing 
$\imath_i^{-1}(B)$ then we attach a disk to $F$ along $B$ and we extend 
$\imath_i$ over $D_i$ for those $i=1,2$ which satisfy the above condition.
By performing these operations for all components of $\p F$ disjoint from 
$\Gamma,$ we extend $O_i$ to an overlap $\bar O_i$ of
$(\Gamma_1,F_1)$ and $(\Gamma_2,F_2)$.
Notice that every overlap of these graphs which restricts to $O_i$ 
must factor through $\bar O_i$. Therefore we proved
 
\begin{corollary}
${\bar O_1},...,{\bar O_d}$ is a basis of overlaps of
$(\Gamma_1,F_1)$ and $(\Gamma_2,F_2)$.
\end{corollary}

%
\section{$A_1$-webs}
\label{s_a1}
%
Interesting examples of confluent and terminal reduction rules 
come from Kuperberg's spider webs associated with simple Lie algebras of rank 
$\leq 2$. These are spaces of graphs in $D^2$ considered modulo certain 
relations, of the type defined in Introduction. 
The classes of graphs considered and the relations between them
depend on the Lie algebra in question.
Because of their relations to quantum invariants, it is important to
consider Kuperberg's webs in surfaces other than 
$D^2,$ even though they are not spiders anymore, since the join operation is 
no longer defined. For that reason Kuperberg's graphs in surfaces 
other than $D^2$ will simply be called {\em webs.}
Although Kuperberg's original reduction rules are not confluent for webs
we will show that these rules can be extended to finite, confluent, 
and terminal sets of reduction rules.

The $A_1$-webs without external vertices are unoriented link diagrams. 
To put such diagrams in the
framework of surface graphs, we define 
{\em crossings} as marked $4$-valent vertices depicted as
$\diag{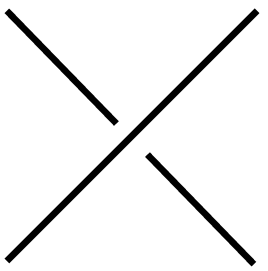}{.3in}$. We require that opposite edges 
meeting at any crossing have equal labels, when taken with coinciding 
orientations.
The notions of linear graphs, reduction rules, 
local and global confluence extend to graphs with crossings in an obvious 
way. Furthermore, the method of
proving local confluence discussed in Section \ref{ss_overlaps} holds for 
graphs with crossings as well.

Consider now a surface $F$ (not necessarily oriented) with a distinguished set
of base points $B\subset \p F$ (possibly empty).
An {\em $A_1$-web} in $(F,B)$ is an unoriented 
graph all of whose internal vertices are crossings
and all of whose external vertices are points of $B$.
(Such graphs in $D^2$ are called {\em unoriented tangle diagrams} with 
endpoints in $B$.) 
We denote the set of all $A_1$-webs in $(F,B)$ by $\W_{A_1}(F,B)$.
Let $R$ be a fixed ring with a distinguished invertible element $A$.
The {\em $A_1$-web space} over $R$ is the $R$-module 
$$\A_1(F,B,R)= R \W_{A_1}(F,B)/\R(T_1,T_2),$$
where

$$\begin{array}{lll}
T_1 & = & \diag{crossa}{.3in}- A\diag{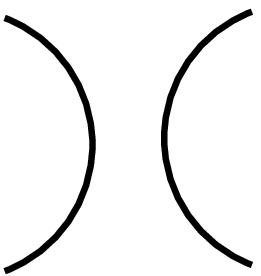}{.3in}-A^{-1}
\diag{smoothb}{.3in}\\
T_2 & = & \diag{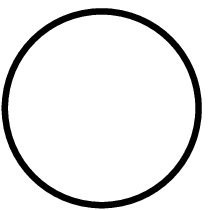}{.3in} -(A^2+A^{-2})\emptyset.
\end{array}$$

Here and further on, all relations take place in $D^2$ and all $1$-valent 
vertices are external, unless stated otherwise.
The above relations suggest the obvious reduction rules:
$$S_1:\diag{crossa}{.3in}\to A\diag{smootha}{.3in}+A^{-1}
\diag{smoothb}{.3in} \quad,\hspace*{.3in} S_2:\diag{circle}{.3in}\to 
-(A^2+A^{-2})\emptyset.$$
Denote the number of crossings and connected components of 
$\Gamma\in \W_{A_1}(F,B)$ by $v(\Gamma)$ and $c(\Gamma),$ respectively.
If $\Z_{\geq 0}\times \Z_{\geq 0}$ is given the lexicographic ordering then
these reduction rules replace each graph $\Gamma$ by a combination of graphs
$\Gamma_i$ such that $(v(\Gamma_i),c(\Gamma_i))<(v(\Gamma),c(\Gamma))$.
Therefore, the rules $S_1,S_2$ are terminal. The irreducible graphs are 
those with no crossings and no contractible loops.
Since \diag{crossa}{.3in} and \diag{circle}{.3in} have no non-trivial 
overlaps, they are locally confluent and, hence, also globally confluent.
Now Theorem \ref{main}(2) provides answers to questions (1),(2)
of Introduction:

\begin{corollary}\label{A_1basis}
For any $F,R$ and $q,$ 
$\A_1(F,B,R)$ is the free $R$-module with a basis given by finite collections
of disjoint non-trivial simple closed loops in $F,$ including $\emptyset$.
\end{corollary}

%
\section{$A_2$-webs}
\label{s_a2}
%
Let $F$ be a surface with a distinguished set of 
base points $B \subset \p F$ (possibly empty) 
which are marked by $\pm 1$.
An {\em $A_2$-web} in $(F,B)$ is an oriented graph $\Gamma$ in $F$ 
all of whose internal vertices are either $3$-valent 
sinks or sources or $4$-valent crossings:

\begin{center}
\diag{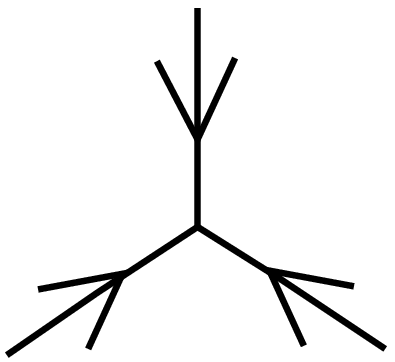}{.4in}\hspace*{.4in} \diag{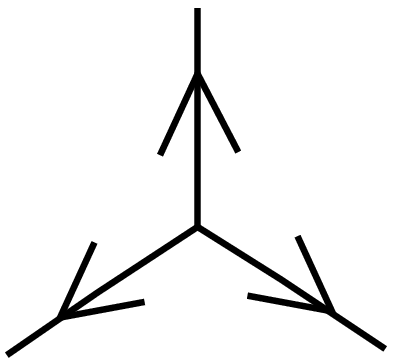}{.4in}\hspace*{.4in}
\diag{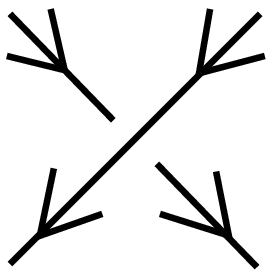}{.3in}\hspace*{.4in} \diag{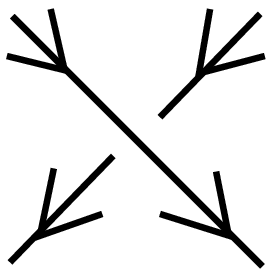}{.3in}
\end{center}
and such that all external vertices are points of $B$.
Furthermore, we require that the external edge adjacent to $b\in B$
is oriented inwards or outwards according to the 
labeling of $b$ by $1$ or $-1$.

Denote the set of $A_2$-webs in $F$ by $\W_{A_2}(F,B)$. 
The {\em $A_2$-web space} is 
$$\A_2(F,B,R)= R \W_{A_2}(F,B)/\R(T_1,T_2,T_3,T_4,T_5,T_6),$$
where
\renewcommand{\arraystretch}{2}
$$\begin{array}{lll}
T_1 & = & \diag{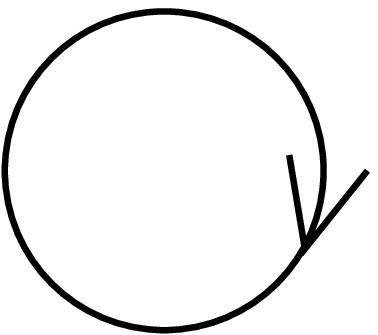}{.3in}\ -\ (q+1+q^{-1})\emptyset,\\
T_2 & = & \diag{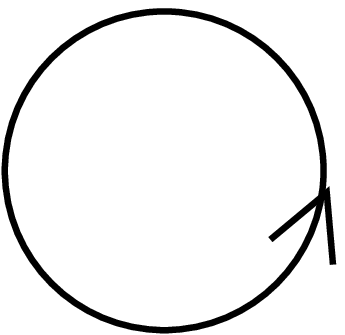}{.3in}\ -\ (q+1+q^{-1})\emptyset,\\
T_3 & = & \parbox{.8in}{\psfig{figure=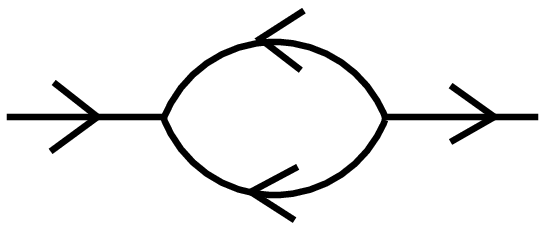,height=.3in}}+
    (q^\frac{1}{2}+q^{-\frac{1}{2}})
\parbox{.6in}{\psfig{figure=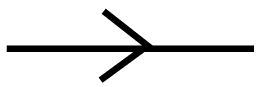,height=.1in}}\\
T_4 & = &
\diag{square}{.45in}-\diag{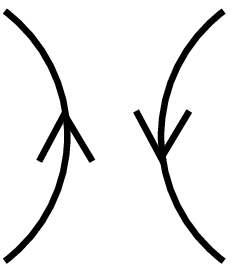}{.35in}-\diag{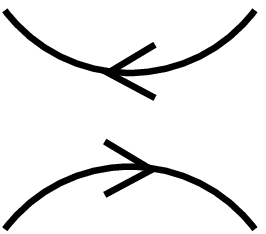}{.35in}\\
T_5 & = & \diag{crossp}{.3in}-  q^{\frac{1}{6}}
\parbox{.4in}{\psfig{figure=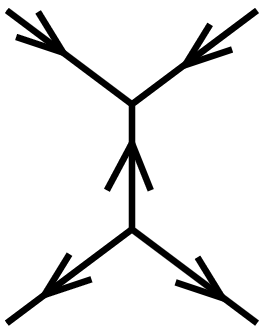,height=.4in}}-
q^{-\frac{1}{3}}\diag{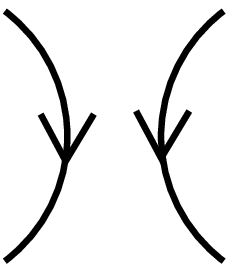}{.3in}\\
T_6 & = & \diag{crossm}{.3in} -  q^{-\frac{1}{6}}
\parbox{.4in}{\psfig{figure=crossh,height=.4in}}-
q^{\frac{1}{3}}\diag{smootho}{.3in},\\
\end{array}$$
and $R$ is an arbitrary ring with a distinguished invertible element 
denoted by $q^{\pm \frac{1}{6}}\in R$. ($R=\C[q^{\pm \frac{1}{6}}]$ in 
\cite{Ku-spider}.)

The above relations suggest ``obvious'' reduction rules:

$$S_1:\ \diag{circleclock}{.3in}\ \to (q+1+q^{-1})\emptyset,
\quad S_2:\ \diag{circlecounter}{.3in} \to
(q+1+q^{-1})\emptyset,$$
$$S_3:\ \parbox{.8in}{\psfig{figure=bigon,height=.3in}}\ \to
-(q^\frac{1}{2}+q^{-\frac{1}{2}})
\parbox{.6in}{\psfig{figure=line,height=.1in}}$$
$$S_4:\ \diag{square}{.5in}\ \to
\diag{sqsmooth1}{.4in}+\diag{sqsmooth2}{.4in}$$
$$S_5:\ \diag{crossp}{.3in}\to  q^{\frac{1}{6}}
\parbox{.4in}{\psfig{figure=crossh,height=.4in}}+
q^{-\frac{1}{3}}\diag{smootho}{.3in}\quad ,\quad
S_6:\ \diag{crossm}{.3in} \to q^{-\frac{1}{6}}
\parbox{.4in}{\psfig{figure=crossh,height=.4in}}+
q^{\frac{1}{3}}\diag{smootho}{.3in}.$$

Denote the number of connected components, $3$-valent vertices, and
crossings of any $\Gamma \in \W_{A_2}(F,B)$ by $c(\Gamma),$
$v_3(\Gamma),$ and $v_4(\Gamma),$ respectively.
If $\Z_{\geq 0}\times \Z_{\geq 0}\times \Z_{\geq 0}$ is given the
lexicographic ordering, then the above reduction rules 
replace $\Gamma$ by a linear combination of graphs $\Gamma_i$ such that
$(v_4(\Gamma_i),v_3(\Gamma_i),c(\Gamma_i))<(v_4(\Gamma),v_3(\Gamma),
c(\Gamma))$. Consequently, these reduction rules are terminal.
However, they are not confluent for $F\ne D^2, S^2$! 
Indeed, any surface $F\ne D^2,S^2$ contains
an annulus whose core is not contractible in $F$ and the two possible
applications of $S_4$ 
to $\Gamma=\diag{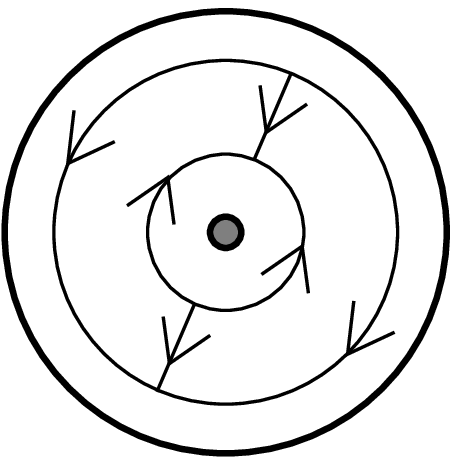}{.7in},$ followed by $S_1, S_2,$
reduce it to
$X_1+(q+1+q^{-1})\emptyset$ and $X_2+(q+1+q^{-1})\emptyset,$ where
$$X_1=\diag{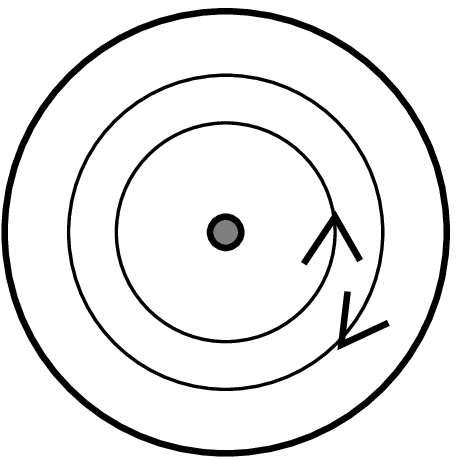}{.6in},\quad X_2=\diag{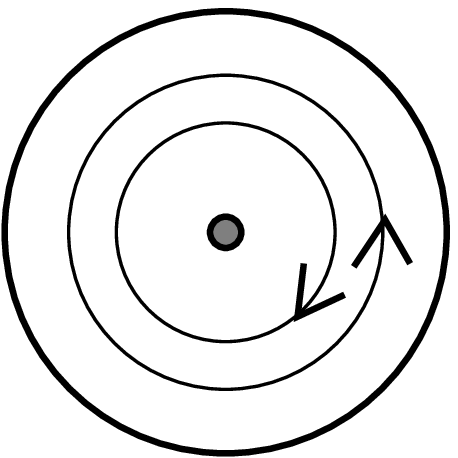}{.6in}.$$
Since $X_1$ and $X_2$  are irreducible and not isotopic,
the reduction rules are not confluent.
In order to remedy this imperfection, we need to consider an additional
reduction rule taking place in an annulus:
$$S_7:\ \diag{weakiso1}{.6in}\to \diag{weakiso2}{.6in}.$$

Note that $S_1,...,S_7$ are terminal in the oriented case but
not terminal in the unoriented case, c.f. last paragraph of Section
\ref{ss_mfld_graphs}. Indeed, $S_7$ is its own inverse
in the unoriented case! Therefore, we consider
$A_2$-webs in oriented surfaces only and work in the category of oriented 
surfaces from now on. (We assume that all surfaces appearing in
reduction rules $S_1,...,S_7$ have counterclockwise orientation.)

One checks all overlaps for $S_1,...,S_7$ and
concludes that $S_1,...,S_7$ are locally confluent on all of them\footnote{
Recall that irreducible overlaps of $S_3$ and $S_4$ (as abstract graphs) were
classified in Example \ref{overlap_example}.}.
Therefore, by Theorem \ref{main}, we conclude:

\begin{corollary}\label{a2_main}
The reduction rules $S_1,...,S_7$ are both terminal and confluent for
graphs in $\W_{A_2}(F,B),$ for any oriented surface $F$ and any set of 
marked base 
points $B\subset \p F$. Consequently, $\A_2(F,B,R)$ is a free $R$-module 
with a basis composed of irreducible graphs in $\W_{A_2}(F,B)$.
\end{corollary}

Observe that irreducible
$A_2$-webs in $D^2$ are those which have 
no $S^1$'s, no internal bi-gons, and no internal $4$-gons.
(Such graphs in $D^2$ are called {\em non-elliptic}
in \cite{Ku-spider}.)
While these terms are intuitively obvious for graphs in $D^2,$ they do
require clarification for graphs in other surfaces.

Components of $F\setminus \Gamma$ are {\em faces of $\Gamma$}.
A face is {\em internal} if it is disjoint from $\p F$. 
An internal face is called an {\em $n$-gon} if it is a disk bounded by  
a sequence of $n$ edges of $\Gamma$. (The orientations of the 
edges are irrelevant.) An $n$-gon is {\em true} if 
all its boundary edges are distinct; otherwise it is {\em fake}.
For example, the $4$-gon in $S^1\times I$ bounded by the edges
$E_1,E_2,E_3,E_2$ depicted below is fake.

$$\parbox{1in}{\psfig{figure=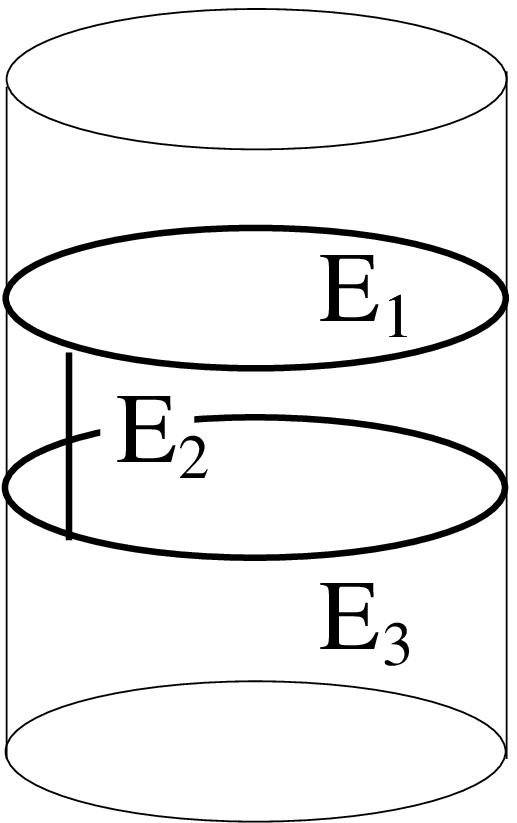,height=1.3in}}$$

A loop in $\Gamma$ bounding a disk in $F\setminus \Gamma$ 
is called a $0$-gon. The next statement follows directly from $A_2$-web 
reduction rules:

\begin{corollary}\label{A_2basis}
The irreducible graphs in $\W_{A_2}(F,B)$ are precisely those
with no $0$-gons, no true bi-gons, and no true $4$-gons.
\end{corollary}

%
\section{$B_2$-webs}
\label{s_b2} 
%

Let $F$ be a surface together with a specified finite set of base points
$B \subset \p F$ (possibly empty), each of them marked by 
$1$ or $2$. Throughout this section we work in the category of unoriented
surfaces.

\begin{definition}\label{b2-web_def} 
Let $\W_{B_2}(F,B)$ be the set of all labeled graphs $\Gamma$ in $F$,
with $\Lambda=\{1,2\}$ and $\nu$ the identity, c.f. Section \ref{s_graphs},
such that\\
(1) the labels of edges adjacent to points of $B$ coincide with their labels,
and\\
(2) all internal vertices of $\Gamma$ are of the form
$$\diag{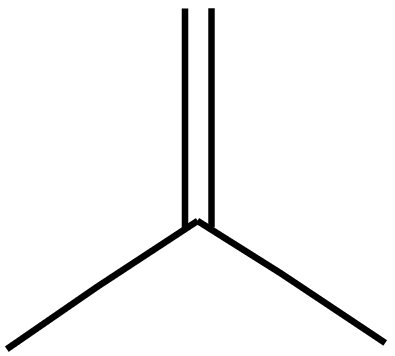}{.4in}$$
(The labels are depicted either by a single or double line.)
\end{definition}

The {\em $B_2$-web space} is
$$\B_2(F,B,R)=R\W_{B_2}(F,B)/\R(T_1,...,T_6),$$
where

\renewcommand{\arraystretch}{2}
$$\begin{array}{lll}
T_1 & = & \diag{circle}{.3in}+(q^2+q+q^{-1}+q^{-2})\emptyset\\
T_2 & = & \diag{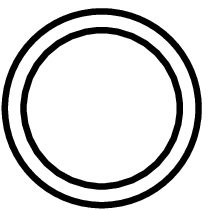}{.3in}-(q^3+q+1+q^{-1}+q^{-3})\emptyset\\
T_3 & = & \parbox{.8in}{\psfig{figure=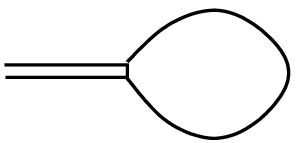,height=0.25in}}\\
T_4 & = & \parbox{.8in}{\psfig{figure=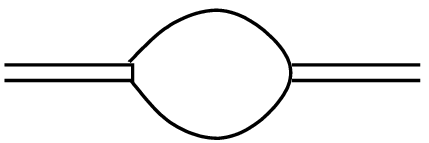,height=0.25in}}+(q+2+q^{-1})
    \parbox{1in}{\psfig{figure=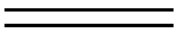,height=0.05in}}\\
T_5 & = & \parbox{.6in}{\psfig{figure=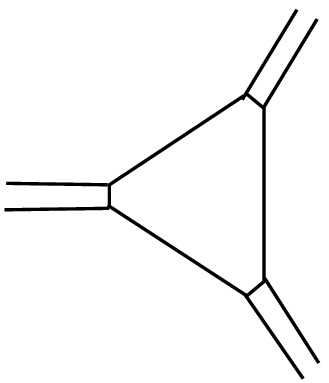,height=0.5in}}\\
T_6 & = & \parbox{.4in}{\psfig{figure=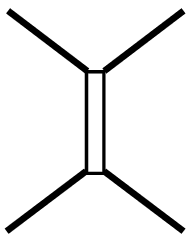,height=0.35in}}-
    \parbox{.4in}{\psfig{figure=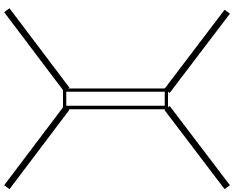,height=0.27in}}+\diag{smootha}{.3in}-
   \parbox{.3in}{\psfig{figure=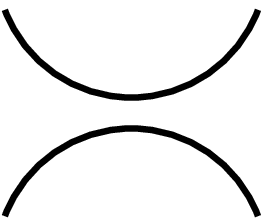,height=0.27in}}
\end{array}$$
and $R$ is an arbitrary ring with a distinguished invertible element
$q\in R$.
$B_2$-webs with crossings are discussed in the next
subsection.

While $T_1,...,T_5$ yield ``obvious'' reduction relations, which we
denote by $S_1,...,S_5,$ relation $T_6$ creates a problem since
the rule
$$\parbox{.4in}{\psfig{figure=b2e,height=0.3in}}\to
    \parbox{.4in}{\psfig{figure=b2f,height=0.27in}}-\diag{smootha}{.3in}+
   \parbox{.3in}{\psfig{figure=smoothb,height=0.27in}}$$
is its own inverse and, hence, it is not terminal.
Following Kuperberg's idea, we remedy this problem by allowing
$B_2$-webs to have $4$-valent vertices subject to a relation $T_6'=0,$ where
$$T_6'=\parbox{.4in}{\psfig{figure=b2e,height=0.3in}}-
\diag{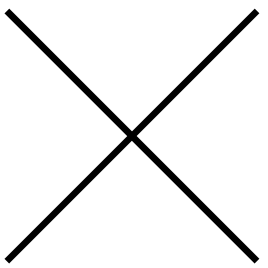}{.3in}- \diag{smoothb}{.3in}\quad .$$
We denote this extended family of webs by $\W_{B_2}'(F,B)$. 
Note that $T_6'$ does not introduce any new relations and, therefore,
 $$\B_2(F,B,R)=R\W_{B_2}(F,B)/\R(T_1,...,T_6)=
R\W_{B_2}'(F,B)/\R(T_1,...,T_5,T_6').$$ 
Now, $T_6'$ suggests the reduction rule
$$S_6: \parbox{.4in}{\psfig{figure=b2e,height=0.3in}}\to 
\diag{cross}{.3in}+ \diag{smoothb}{.3in}$$

Since each of the reduction rules $S_1,...,S_6$ either decreases 
the number of vertices or decreases the number of connected components
without increasing the number of vertices, 
these reduction rules are terminal. However, they are not 
confluent in general\footnote{
These rules may be confluent for
certain choices of $R$ and $q,$ but not for all.}.
$S_1,...,S_6$ have the following basis of overlaps:
$O_{36}=\parbox{0.8in}{\psfig{figure=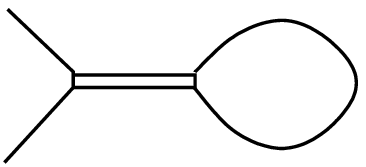,height=0.3in}},\  
O_{46a}=\parbox{1in}{\psfig{figure=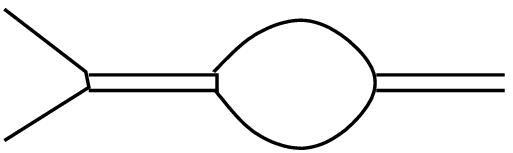,height=0.25in}},\ 
O_{46b}=\parbox{.8in}{\psfig{figure=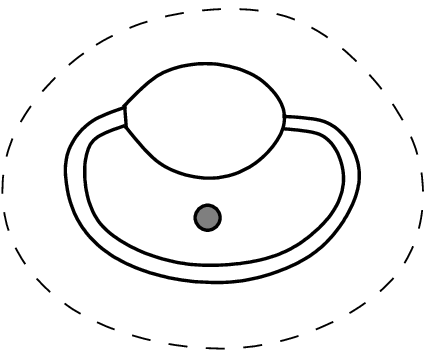,height=0.6in}},\\
O_{46c}=\parbox{.8in}{\psfig{figure=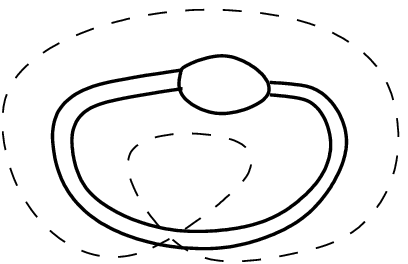,height=0.4in}},\ 
O_{56a}=\parbox{.6in}{\psfig{figure=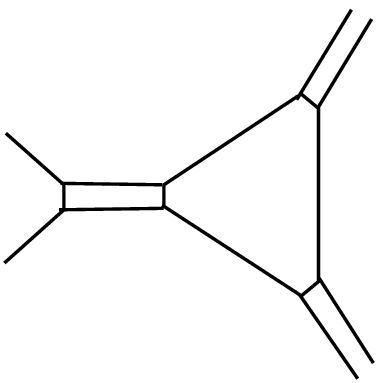,height=0.5in}},
O_{56b}=\parbox{.7in}{\psfig{figure=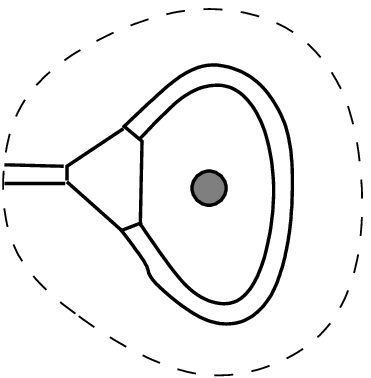,height=0.6in}},
O_{56c}=\parbox{.6in}{\psfig{figure=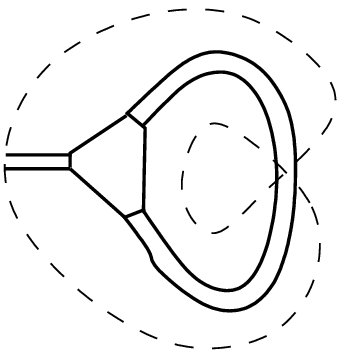,height=0.6in}}.$\\
Dashed lines denote boundaries of surfaces. (No dashed line is drawn for
diagrams in $D^2$.)
Overlaps $O_{46b}$ and $O_{56b}$ take place in annuli and
$O_{46c}$ and $O_{56c}$ in M\"obius bands.
Unfortunately, $S_1,...,S_6$ are not locally confluent on 
these overlaps. For example, 

$$0 \stackrel{S_3}{\longleftarrow}
\parbox{0.8in}{\psfig{figure=b2a2,height=0.3in}}
\stackrel{S_6,S_1}{\longrightarrow} \diag{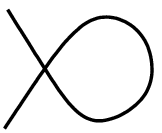}{.25in}\ -
(q^2+q+q^{-1}+q^{-2})\hspace*{.1in} 
\parbox{.2in}{\psfig{figure=linev, height=.3in}},$$
and both of these linear graphs are irreducible with respect to
$S_1,...,S_6$. In order to remedy that, we introduce the
following new rules (preserving relations $S_1$-$S_6$):
$$\begin{array}{llll}
S_7: & \diag{kinku}{.25in} & \to &
(q^2+q+q^{-1}+q^{-2})\hspace*{.1in} 
\parbox{.2in}{\psfig{figure=linev, height=.3in}}\\
S_8: & \parbox{.4in}{\psfig{figure=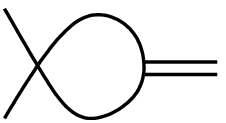, height=.25in}} & \to &
-(q+2+q^{-1})\ \diag{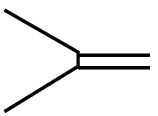}{.25in}\\
S_9: & \parbox{.6in}{\psfig{figure=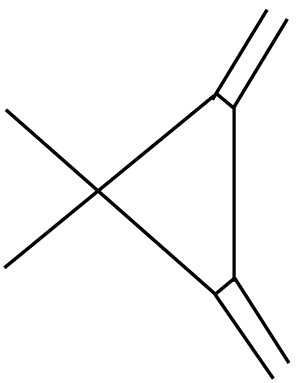, height=.5in}} & \to &
(q+2+q^{-1})\ \diag{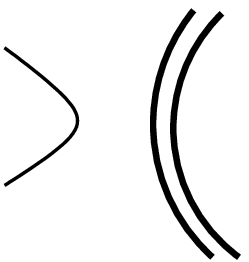}{.3in}\\
S_{10}: & \diag{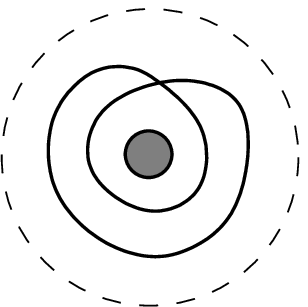}{.4in}  & \to & -(q+2+q^{-1})\ 
\diag{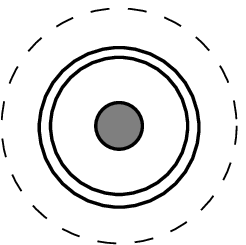}{.4in}+(q^2+q+q^{-1}+q^{-2})\emptyset\\
S_{11}: & \diag{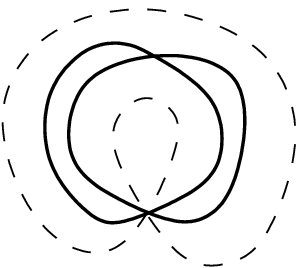}{.45in}  & \to & -(q+2+q^{-1})\ 
\diag{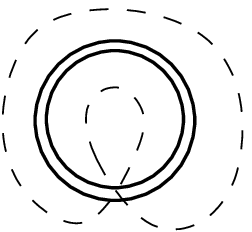}{.45in}+(q^2+q+q^{-1}+q^{-2})\emptyset\\
S_{12}: & \diag{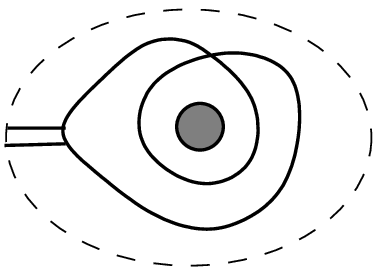}{.5in}  & \to & 0\\
S_{13}: & \diag{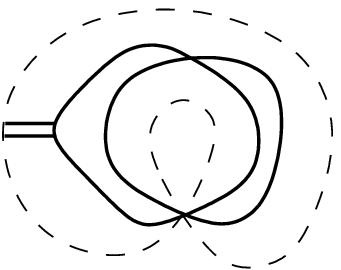}{.5in}  & \to & 0\\
\end{array}$$

Rules $S_{10},S_{12}$ take place in annuli and rules
$S_{11},S_{13}$ in M\"obius bands. (The graph on the left
side of rule $S_{11}$ has a single vertex, and the one on the
left side of $S_{13}$ has two vertices.)

Now $S_1,...,S_{13}$ are locally confluent on
$O_{36},...,O_{56c},$ but the new rules, $S_7$-$S_{13},$ lead to new 
overlaps:\\
$O_{77}=\parbox{.7in}{\psfig{figure=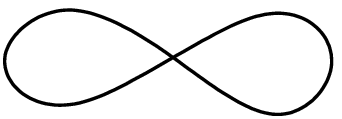,height=0.20in}},\ 
O_{68}=\parbox{.6in}{\psfig{figure=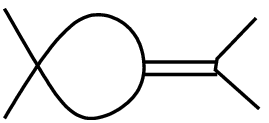,height=0.25in}},\ 
O_{78}=\parbox{.6in}{\psfig{figure=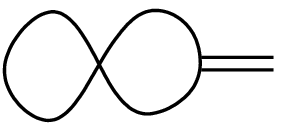,height=0.25in}},\ 
O_{88}=\parbox{.8in}{\psfig{figure=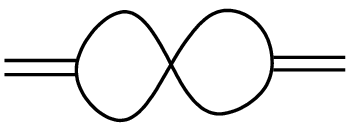,height=0.25in}},\\
O_{69a}=\parbox{.5in}{\psfig{figure=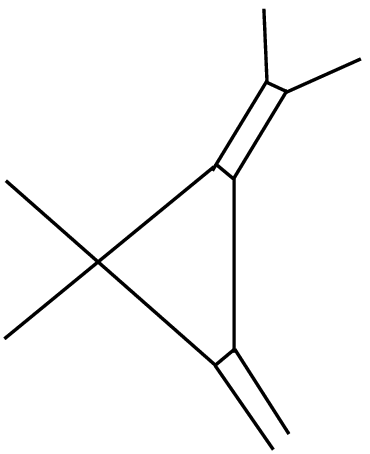,height=0.6in}},\ 
O_{69b}=\parbox{.6in}{\psfig{figure=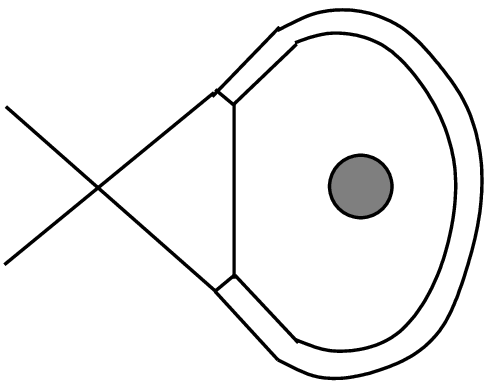,height=0.4in}},\ 
O_{69c}=\parbox{.6in}{\psfig{figure=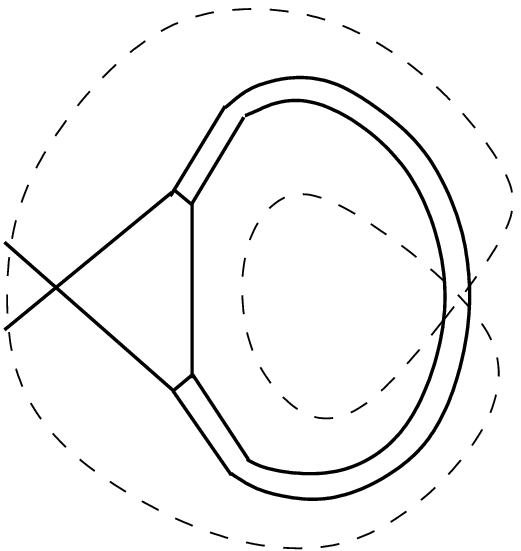,height=0.6in}},\ 
O_{79}=\parbox{.6in}{\psfig{figure=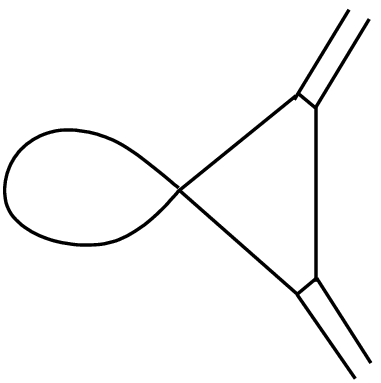,height=0.5in}},\\
O_{89}=\parbox{.7in}{\psfig{figure=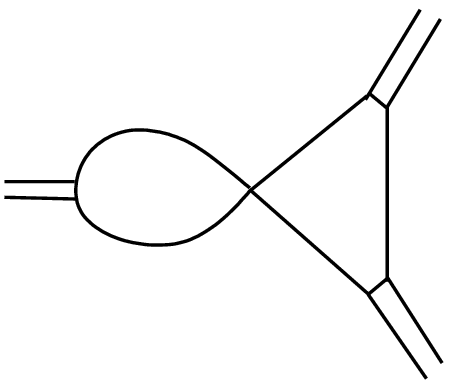,height=0.5in}},\ 
O_{99}=\parbox{.6in}{\psfig{figure=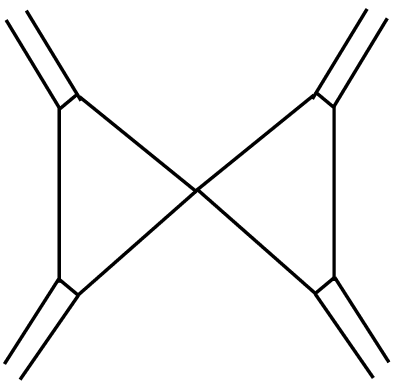,height=0.5in}},\ 
O_{6,12}=\parbox{.8in}{\psfig{figure=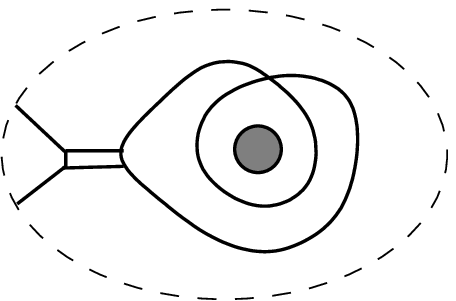,height=0.5in}},\ 
O_{6,13}=\parbox{.8in}{\psfig{figure=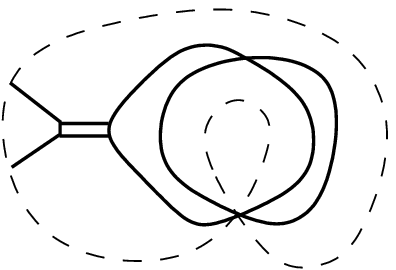,height=0.5in}}$.
Now, we make $S_1,...,S_{13}$ locally confluent on
these overlaps by introducing the following new reduction rules:\\

$S_{14}:\ \parbox{.6in}{\psfig{figure=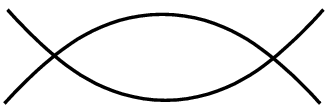,height=0.15in}}\
\to\ -(q+2+q^{-1})\diag{cross}{.3in}-(q^2+2q+2+2q^{-1}+q^{-2})
\diag{smootha}{.3in}\\
S_{15}:\ \parbox{.6in}{\psfig{figure=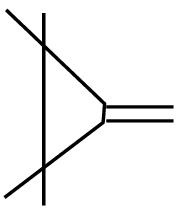,height=0.35in}}\ \to\ 
(q+2+q^{-1})\left(\diag{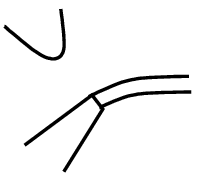}{.35in}+\diag{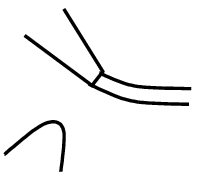}{.35in}\right)\\
S_{16}:\ \parbox{.6in}{\psfig{figure=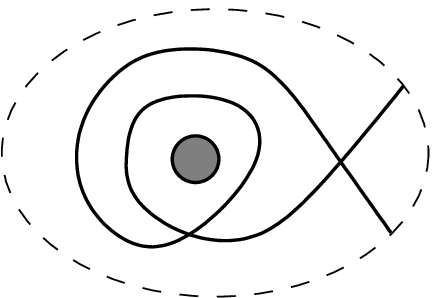,height=0.4in}}\ \to \ (q+2+q^{-1})\ 
\parbox{.6in}{\psfig{figure=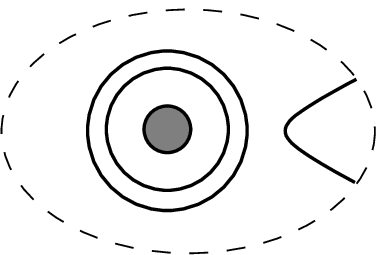, height=.4in}}-
(q^2+q+q^{-1}+q^{-2})\parbox{.7in}{\psfig{figure=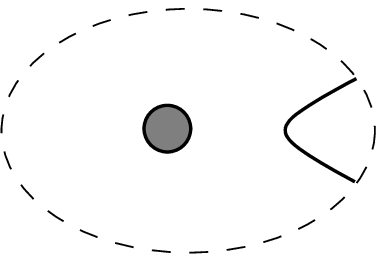, height=.4in}}\\
S_{17}:\ \parbox{.6in}{\psfig{figure=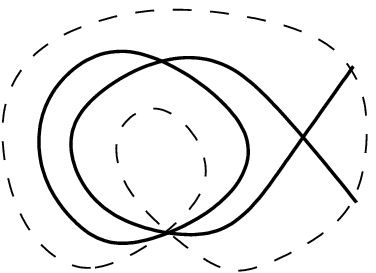,height=0.45in}}\ \to\ 
(q+2+q^{-1})\ \parbox{.6in}{\psfig{figure=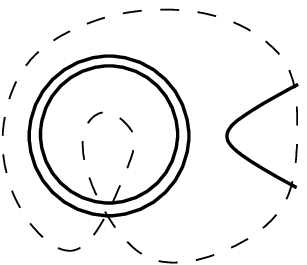, height=.5in}}-
(q^2+q+q^{-1}+q^{-2})\parbox{.6in}{\psfig{figure=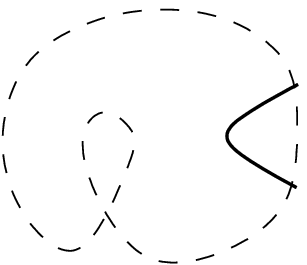, height=.4in}}.\\$

The new overlaps now are:\\
$O_{7,14}=\parbox{.9in}{\psfig{figure=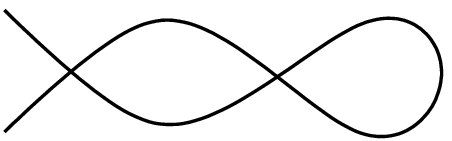,height=0.25in}},\ 
O_{8,14}=\parbox{1in}{\psfig{figure=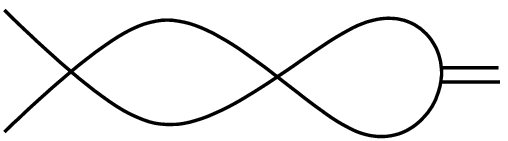,height=0.25in}},\ 
O_{9,14}=\parbox{.7in}{\psfig{figure=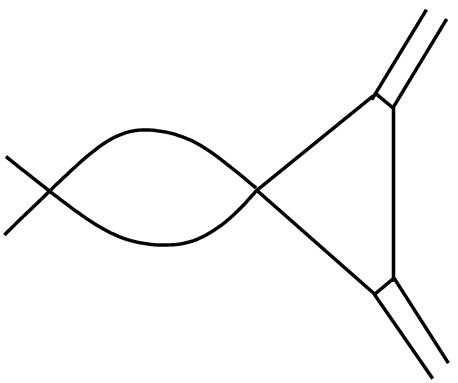,height=0.5in}},\\ 
O_{14,14a}=\parbox{.9in}{\psfig{figure=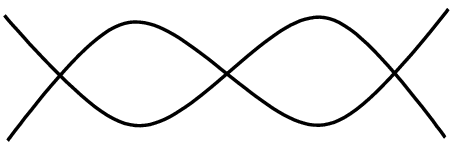,height=0.25in}},\
O_{14,14b}=\parbox{.5in}{\psfig{figure=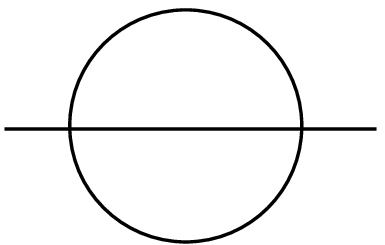,height=0.3in}},\ 
O_{14,14c}=\parbox{.7in}{\psfig{figure=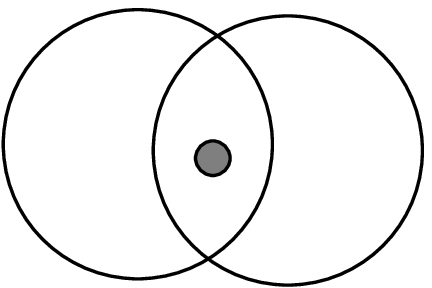,height=0.35in}},\\
O_{14,14d}=\parbox{.6in}{\psfig{figure=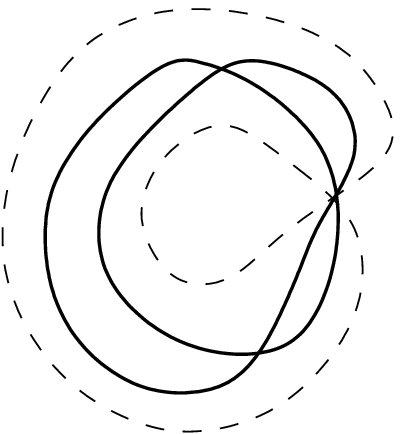,height=0.6in}},\ 
O_{6,15}=\parbox{.6in}{\psfig{figure=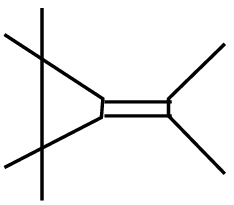,height=0.4in}},\ 
O_{7,15}=\parbox{.6in}{\psfig{figure=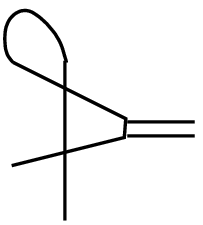,height=0.5in}},\ 
O_{8,15}=\parbox{.6in}{\psfig{figure=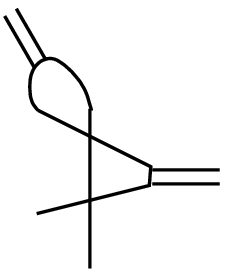,height=0.5in}},\\ 
O_{9,15}=\parbox{.6in}{\psfig{figure=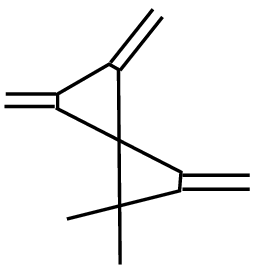,height=0.5in}},\ 
O_{14,15b}=\parbox{.7in}{\psfig{figure=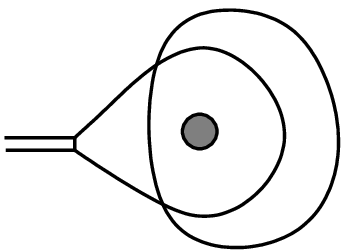,height=0.5in}},\ 
O_{14,15c}=\parbox{.6in}{\psfig{figure=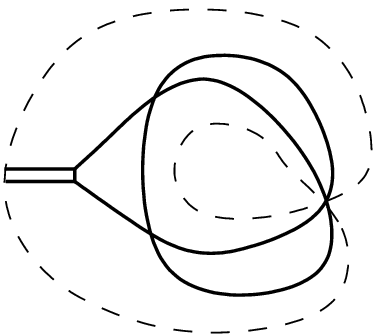,height=0.5in}},\\
O_{15,15a}=\parbox{.7in}{\psfig{figure=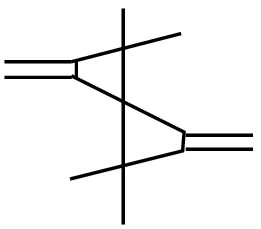,height=0.5in}},\ 
O_{15,15b}=\parbox{.6in}{\psfig{figure=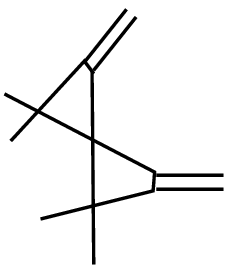,height=0.5in}},\  
O_{15,15c}=\parbox{.8in}{\psfig{figure=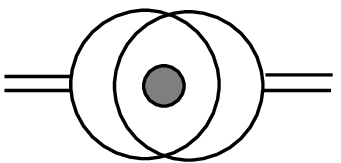,height=0.3in}},\\ 
O_{15,15d}=\parbox{.8in}{\psfig{figure=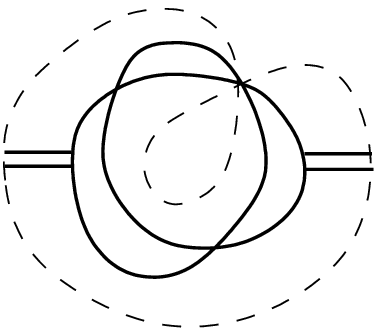,height=0.6in}},\
O_{15,15e}=\parbox{.9in}{\psfig{figure=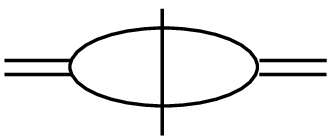,height=0.3in}},\
O_{7,16}=\parbox{.9in}{\psfig{figure=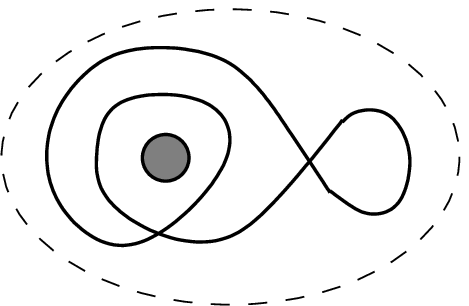,height=0.5in}},\\ 
O_{8,16}=\parbox{.9in}{\psfig{figure=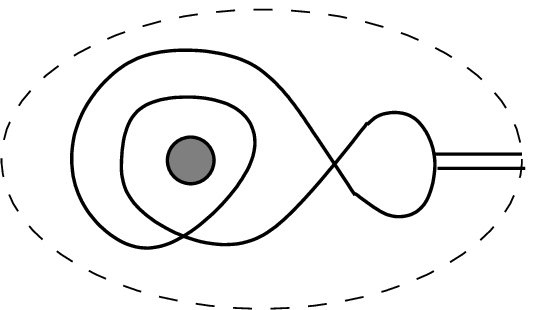,height=0.5in}},\ 
O_{9,16}=\parbox{1in}{\psfig{figure=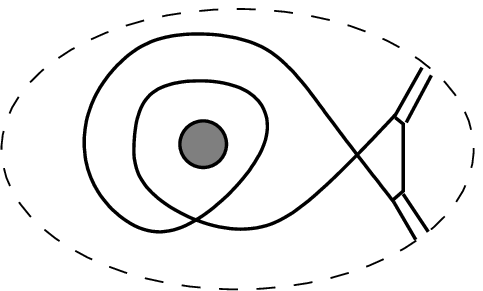,height=0.5in}},\  
O_{15,16}=\parbox{.9in}{\psfig{figure=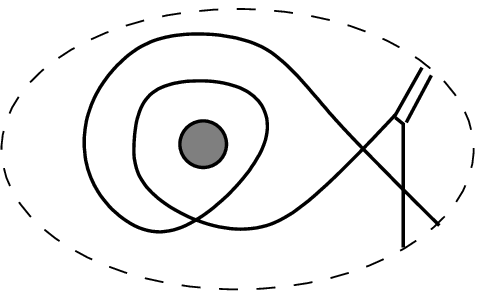,height=0.5in}},\\ 
O_{16,16a}=\parbox{1in}{\psfig{figure=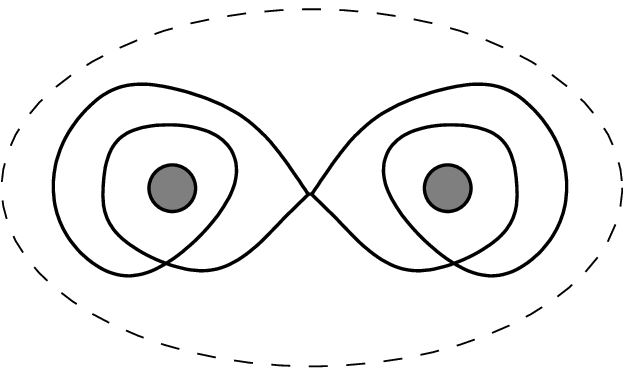,height=0.5in}},\ 
O_{16,16b}=\parbox{.8in}{\psfig{figure=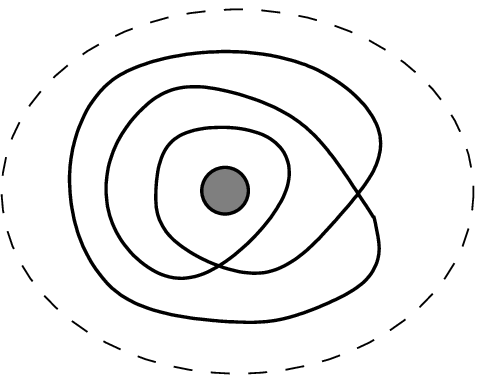,height=0.5in}},\
O_{7,17}=\parbox{.8in}{\psfig{figure=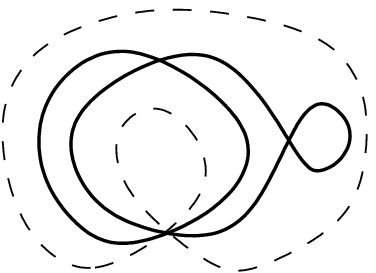,height=0.5in}},\\
O_{8,17}=\parbox{.8in}{\psfig{figure=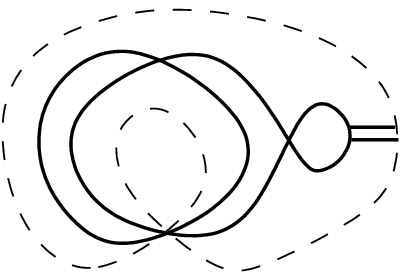,height=0.5in}},\ 
O_{9,17}=\parbox{.9in}{\psfig{figure=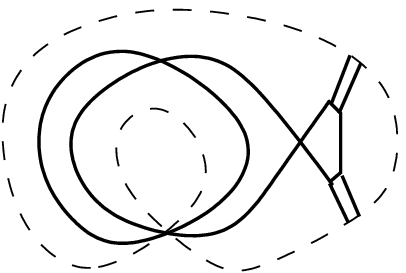,height=0.5in}},\  
O_{14,17}=\parbox{.9in}{\psfig{figure=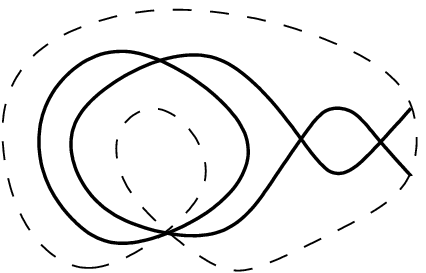,height=0.5in}},\\ 
O_{15,17}=\parbox{.9in}{\psfig{figure=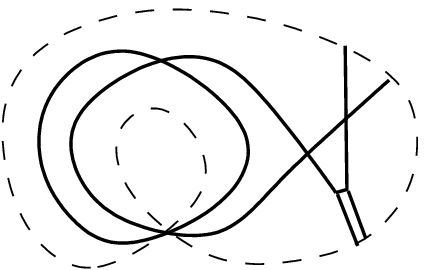,height=0.5in}},\ 
O_{16,17a}=\parbox{1.1in}{\psfig{figure=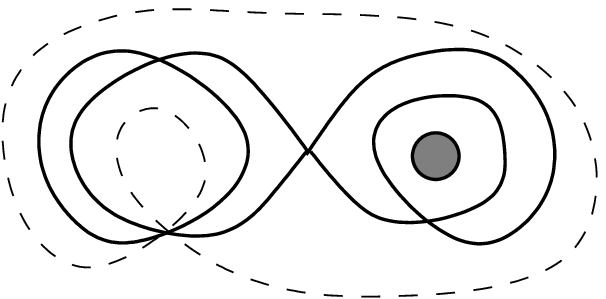,height=0.5in}},\  
O_{16,17b}=\parbox{1in}{\psfig{figure=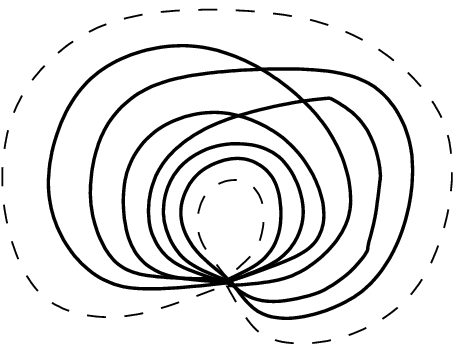,height=0.7in}},\\
O_{17,17a}=\parbox{1.3in}{\psfig{figure=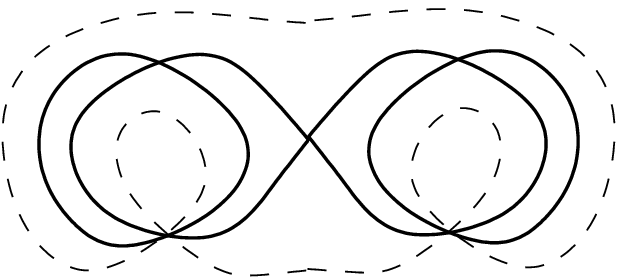,height=0.5in}}$.
Finally, there is an overlap $O_{17,17b}$ in a Klein bottle obtained
by gluing two M\"obius bands along their boundaries and taking a union
$\Gamma$ of two graphs \parbox{.7in}{\psfig{figure=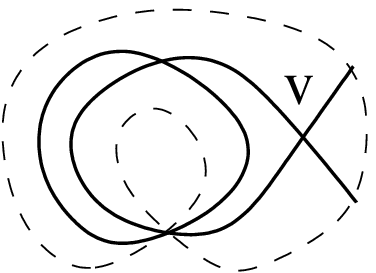,height=0.45in}}
overlapping at vertex $v$. 
(Hence, $\Gamma$ has three $4$-valent vertices and no vertices of
other valences.)
Since $S_1$-$S_{17}$ are not locally confluent on
$O_{6,15},$ we add yet another reduction rule:
$$\begin{array}{llll}
S_{18}:& \parbox{.5in}{\psfig{figure=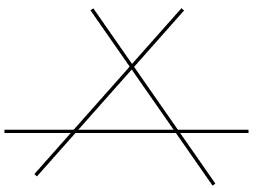,height=.4in}} & \to &
(q+2+q^{-1})\left(\diag{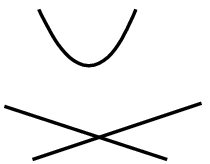}{.4in}\ \ +\diag{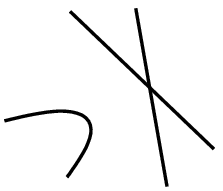}{.4in}\ \
+ \reflectbox{\diag{triangle3}{.4in}} \right)+\\
 & & & (q^2+4q+6+4q^{-1}+q^{-2})\diag{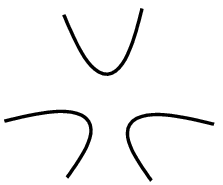}{.4in}\quad .\\
\end{array}$$
One can check that $S_1$-$S_{18}$ are locally confluent on all overlaps listed 
so far. However, $S_{18}$ leads to new overlaps:\\
$O_{7,18}=\parbox{.5in}{\psfig{figure=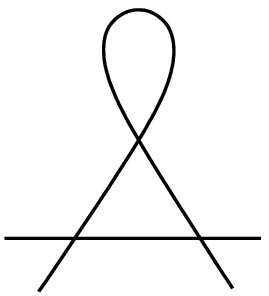,height=0.5in}},\ 
O_{8,18}=\parbox{.5in}{\psfig{figure=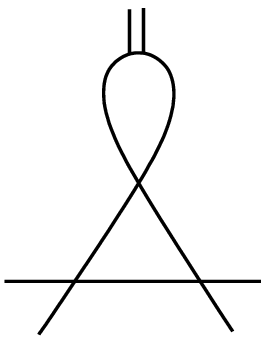,height=0.6in}},\ 
O_{9,18}=\parbox{.5in}{\psfig{figure=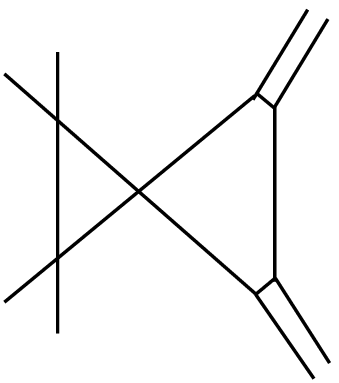,height=0.6in}},\ 
O_{14,18a}=\parbox{.5in}{\psfig{figure=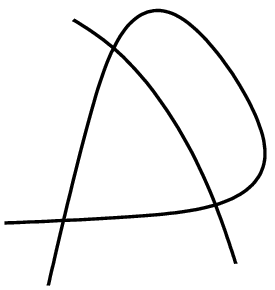,height=0.4in}},\\
O_{14,18b}=\parbox{.7in}{\psfig{figure=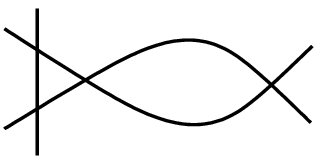,height=0.3in}},\ 
O_{14,18c}=\parbox{.7in}{\psfig{figure=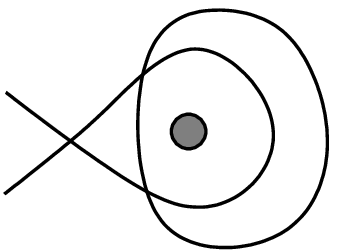,height=0.4in}},\ 
O_{14,18d}=\parbox{.8in}{\psfig{figure=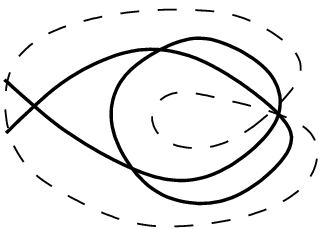,height=0.5in}},\\
O_{15,18a}=\parbox{.5in}{\psfig{figure=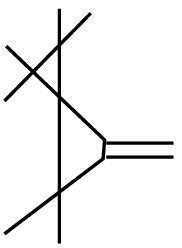,height=0.5in}},\ 
O_{15,18b}=\parbox{.6in}{\psfig{figure=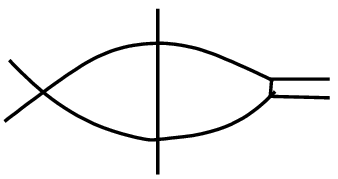,height=0.3in}},\
O_{15,18c}=\parbox{.8in}{\psfig{figure=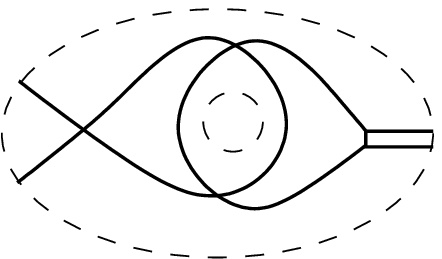,height=0.4in}},\\
O_{15,18d}=\parbox{.8in}{\psfig{figure=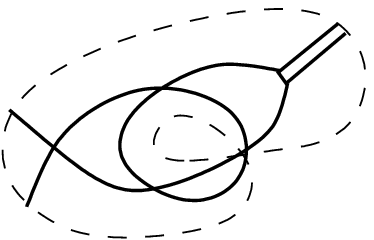,height=0.5in}},\  
O_{15,18e}=\parbox{.9in}{\psfig{figure=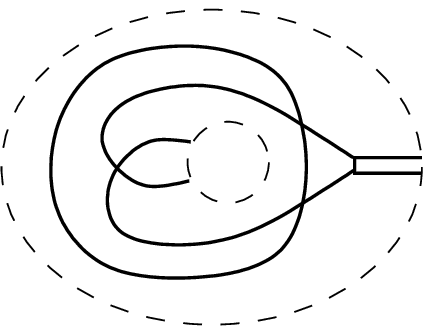,height=0.6in}},\
O_{18,18a}= \parbox{.8in}{\psfig{figure=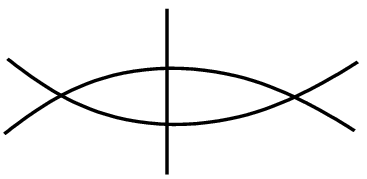,height=0.3in}},\\ 
O_{18,18b}= \parbox{.6in}{\psfig{figure=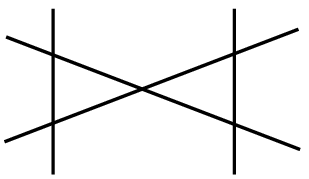,height=0.3in}},\ 
O_{18,18c}= \parbox{.6in}{\psfig{figure=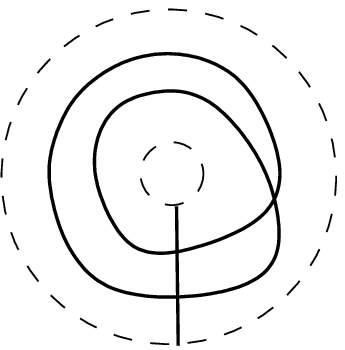,height=0.5in}},\ 
O_{18,18d}= \parbox{.7in}{\psfig{figure=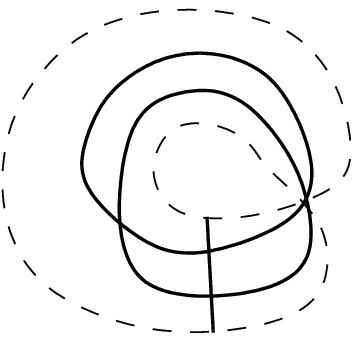,height=0.5in}},\\
O_{18,18d}=\parbox{.8in}{\psfig{figure=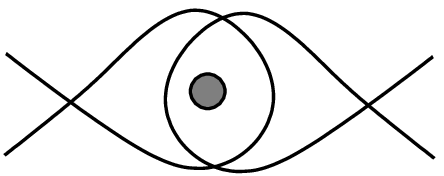,height=0.3in}},\ 
O_{18,18d}= \parbox{.5in}{\psfig{figure=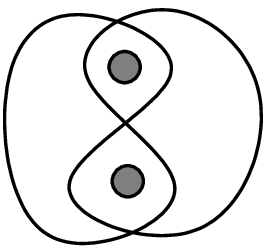,height=0.4in}}.\\$

By checking all of the above overlaps, we conclude that $S_1$-$S_{18}$
are locally confluent on all of them! Therefore, we proved:

\begin{theorem}
The reduction rules $S_1,...,S_{18}$ are both terminal
and confluent for $\W_{B_2}(F,B),$ for any surface $F$ and $B\subset \p F$
(both in orientable and unorientable categories).
\end{theorem}

For sets $B\subset \p F$ all of whose points are marked by $1$,
denote the set of all graphs in $\W_{B_2}(F,B)$ with single edges
only (i.e. edges labeled by $1$), by $\W_{B_2}'(F,B)$. (Graphs in
$\W_{B_2}'(F,B)$ may include double loops.)

\begin{proposition}
For any $F$ and $B\subset \p F$ as above,\\
(1) the embedding $\W_{B_2}'(F,B)\hookrightarrow \W_{B_2}(F,B)$
induces an isomorphism
$$\phi: R\W_{B_2}'(F,B)/\R(S_1, S_2, S_7, S_{10}, S_{11},
S_{14}, S_{16}, S_{17})\to \B_2(F,B,R).$$
(2) the rules $S_1, S_2, S_7, S_{10}, S_{11}, S_{14}, S_{16}, S_{17}$
are terminal and confluent for graphs in $\W_{B_2}'(F,B)$.
\end{proposition}

\begin{proof}
Since internal double edges are resolvable by $S_6,$ $\phi$ is onto.
The map $\phi$ is $1$-$1$ as well:
If $\phi(x)=0$ in $\W_{B_2}(F,B),$ then by confluence of $S_1$-$S_{18}$, 
$x$ can be reduced to $0$ by these rules. Since $x\in R\W_{B_2}'(F,B),$
the rules which do not contain double edges are sufficient to reduce $x$ 
to $0$.
\end{proof}

%
\subsection{$B_2$-webs with crossings}
\label{s_b2cross}
%

For any $F$ and $B\subset \p F$ as in Section \ref{s_b2},
let $\W^c_{B_2}(F,B)$ be the set of all labeled graphs $\Gamma$ in $F$,
with edges labeled by $1$ and $2$ such that\\
(1) the labels of edges adjacent to points of $B$ coincide with their labels,\\
(2) all internal vertices of $\Gamma$ are of the form
$$\diag{b2vertex}{.4in}\ ,\ \diag{crossa}{.3in}\ , \diag{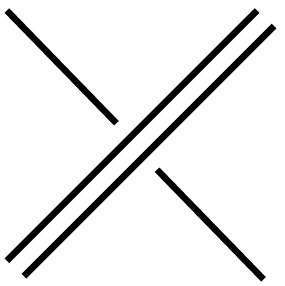}{.3in}\ ,\ 
\diag{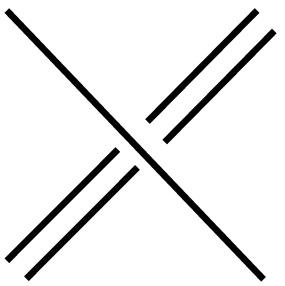}{.3in}\ ,\ \diag{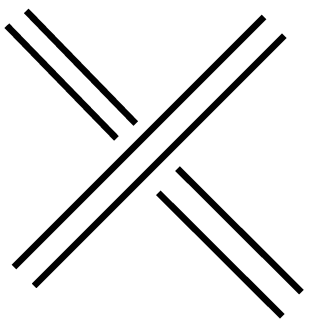}{.3in}.$$
Resolutions of crossings are provided by
$$\begin{array}{llll}
C_1: & \diag{crossa}{.3in} & \to & -q^\frac{1}{2}\diag{smootha}{.3in}-
q^{-\frac{1}{2}}\parbox{.3in}{\psfig{figure=smoothb,height=0.27in}}
+\frac{1}{q^\frac{1}{2}+q^{-\frac{1}{2}}}\diag{cross}{.3in}\\
C_2: & \diag{cross12}{.3in} & \to & 
\frac{q^{-\frac{1}{2}}}{q^\frac{1}{2}+q^{-\frac{1}{2}}}
\diag{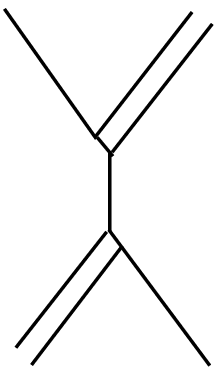}{.5in}+\frac{q^{\frac{1}{2}}}{q^\frac{1}{2}+q^{-\frac{1}{2}}}
\diag{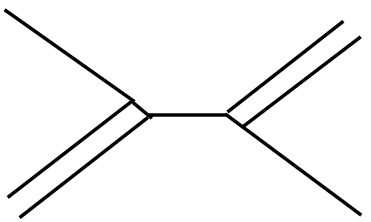}{.3in}\\
C_3: & \diag{cross12b}{.3in} & \to &
\frac{q^{\frac{1}{2}}}{q^\frac{1}{2}
+q^{-\frac{1}{2}}}\diag{fusion1}{.5in}+
\frac{q^{-\frac{1}{2}}}{q^\frac{1}{2}+q^{-\frac{1}{2}}}
\diag{fusion2}{.3in}\\
C_4: & \diag{cross22}{.3in} & \to &
q\diag{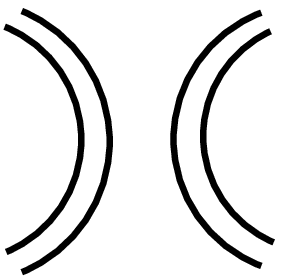}{.3in} + q^{-1}\diag{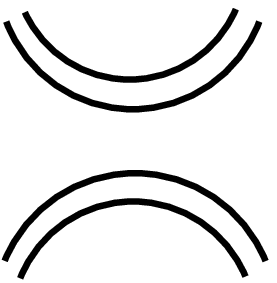}{.3in}
+\frac{1}{q+2+q^{-1}}\diag{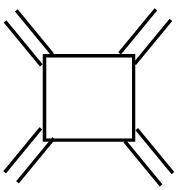}{.35in}.
\end{array}$$
Since the crossing diagrams do not add any new overlaps, we conclude with

\begin{corollary}
If $q^\frac{1}{2}\in R$ is such that $q+1$ is invertible in 
$R,$ then resolutions $S_1,...,S_{18}$ together with $C_1,...,C_4$ 
are confluent and terminal for any surface $F$ and
$B\subset \p F$. 
\end{corollary}

%
\section{$G_2$-webs}
\label{s_g2}
%

Let $F$ be a surface together with a specified finite set of base points
$B \subset \p F$ (possibly empty), each of them marked by 
$1$ or $2$. Throughout this section we work in the category of unoriented
surfaces.

\begin{definition}\label{g2-web_def} 
Let $\W_{G_2}(F,B)$ be the set of all labeled graphs $\Gamma$ in $F$,
whose edges are labeled by $1$ or $2$ and such that\\
(1) the labels of edges adjacent to points of $B$ coincide with their labels,\\
(2) all internal vertices of $\Gamma$ are of the form
$$\diag{b2vertex}{.4in}\qquad or \qquad \diag{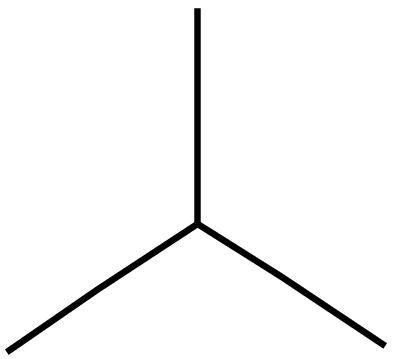}{.4in}.$$
\end{definition}

The {\em $G_2$-web space} is
$${\mathbb G}_2(F,B,R)=R\W_{G_2}(F,B)/\R(S_1,...,S_8),$$
where
\renewcommand{\arraystretch}{2}
$$\begin{array}{lrll}
S_1:& \diag{circle}{.3in} & \to &
(q^5+q^4+q+1+q^{-1}+q^{-4}+q^{-5})\emptyset\\
S_2:& \diag{doublecirc}{.3in} & \to & 
(q^9+q^6+q^5+q^4+q^3+q+2+q^{-1}+q^{-3}+q^{-4}+\\
&  & & q^{-5}+q^{-6}+q^{-9})\emptyset\\
S_3:& \parbox{.6in}{\psfig{figure=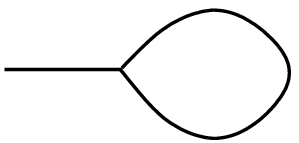,height=0.25in}} & \to & 0\\
S_4:& \parbox{.6in}{\psfig{figure=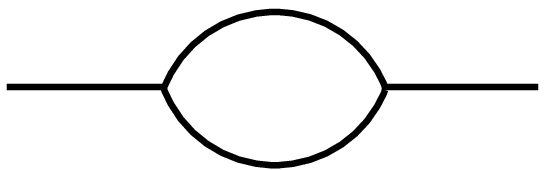,height=0.20in}}& \to & 
-(q^3+q^2+q+q^{-1}+q^{-2}+q^{-3})\ 
\begin{picture}(.3,.1)
\put(0,.05){\line(1,0){.3}}
\end{picture}\\
S_5:& \diag{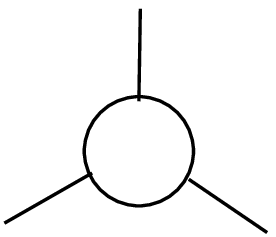}{0.4in} & \to & (q^2+1+q^{-2})\diag{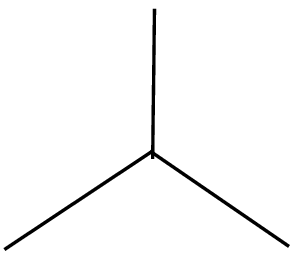}{0.4in}\\
S_6:& \diag{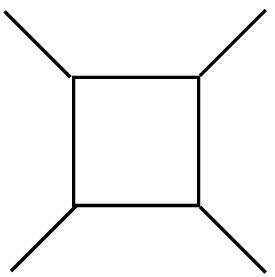}{0.3in} & \to & -(q+q^{-1})
\left(\, \parbox{.3in}{\psfig{figure=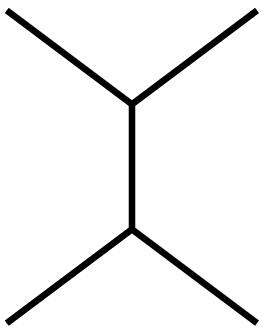,height=.3in}}+
\parbox{.4in}{\begin{sideways} \diag{crosshu}{.3in} 
\end{sideways}}\,  \right)
+(q+1+q^{-1})\left(\diag{smootha}{.3in}+
\parbox{.4in}{\psfig{figure=smoothb,height=.25in}}\right)\\
S_7:&\parbox{.4in}{\psfig{figure=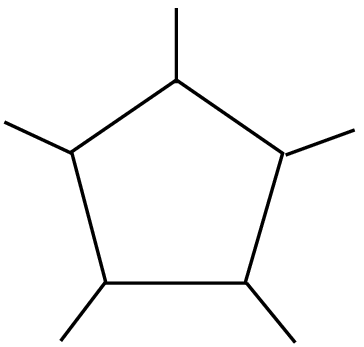,height=0.4in}}& \to &
\left(\parbox{.4in}{\psfig{figure=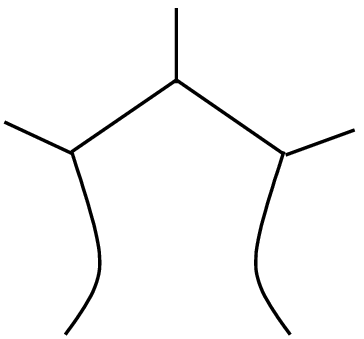,height=0.4in}}+
\parbox{.4in}{\psfig{figure=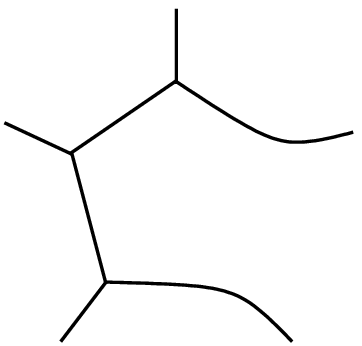,height=0.4in}}+
\parbox{.4in}{\psfig{figure=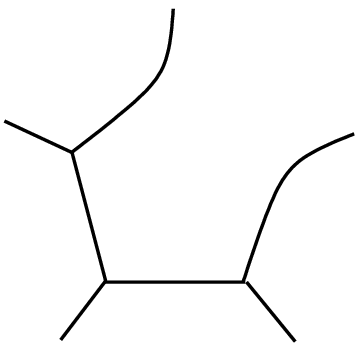,height=0.4in}}+
\reflectbox{\parbox{.4in}{\psfig{figure=pentagon4,height=0.4in}}}+
\reflectbox{\parbox{.4in}{\psfig{figure=pentagon3,height=0.4in}}}
\right)-\\
& & & \left(\parbox{.4in}{\psfig{figure=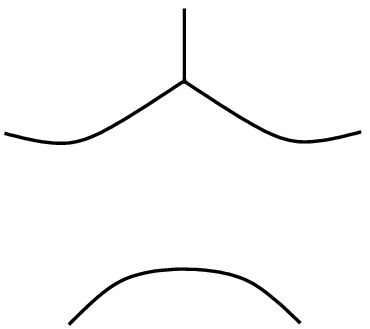,height=0.4in}}+
\parbox{.4in}{\psfig{figure=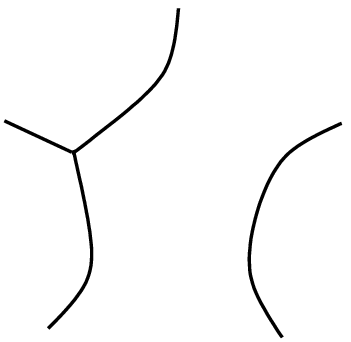,height=0.4in}}+
\parbox{.4in}{\psfig{figure=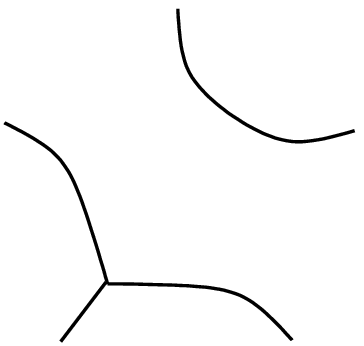,height=0.4in}}+
\reflectbox{\parbox{.4in}{\psfig{figure=penta3,height=0.4in}}}+
\reflectbox{\parbox{.4in}{\psfig{figure=penta2,height=0.4in}}}
\right).\\
S_8:& \parbox{.4in}{\psfig{figure=b2f,height=0.27in}} & \to & 
\parbox{.3in}{\psfig{figure=smoothb,height=0.27in}}-
\parbox{.3in}{\psfig{figure=crosshu,height=.3in}}-
\frac{1}{q^2-1+q^{-2}}\diag{smootha}{.3in}+
\frac{1}{q+1+q^{-1}}\parbox{.4in}{\begin{sideways} \diag{crosshu}{.3in} 
\end{sideways}}.\\
\end{array}$$
(Be advised that reduction rule $S_7$ has wrong signs in
\cite{Ku-spider}.)
Reduction rules for crossings are listed in \cite{Ku-spider}.

By checking all overlaps we conclude

\begin{theorem}
If $F$ is orientable and $q^2-1+q^{-2},q+1+q^{-1}$ are invertible in $R$, then
reduction rules $S_1,...,S_8$ together with
$$S_9: \diag{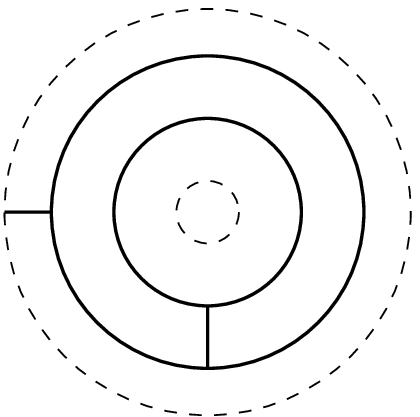}{.6in}\to \diag{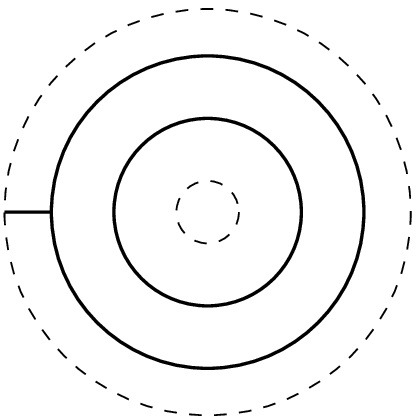}{.6in} -(q+1+q^{-1})
\diag{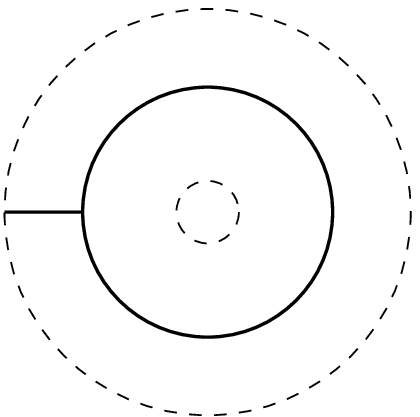}{.5in}$$
(taking place in annulus) are confluent and terminal.\footnote{
We did not check confluence for unoriented surfaces.}
\end{theorem}

%
\section{Partition Category and Dichromatic Reduction Rules}\label{s_dichrom}
%

An example of terminal and confluent reduction rules for abstract graphs
comes from dichromatic polynomial, c.f. \cite{Ye}.

Let $\G_n$ be the set of all unoriented graphs with $n$ external vertices 
labeled from $1$ to $n$ and let 
$R=\Z[p,q,s,v,w_1,w_2]$. Consider reduction rules 
$$\begin{array}{llll}
S_{k,l}:& \parbox{0.8in}{\psfig{figure=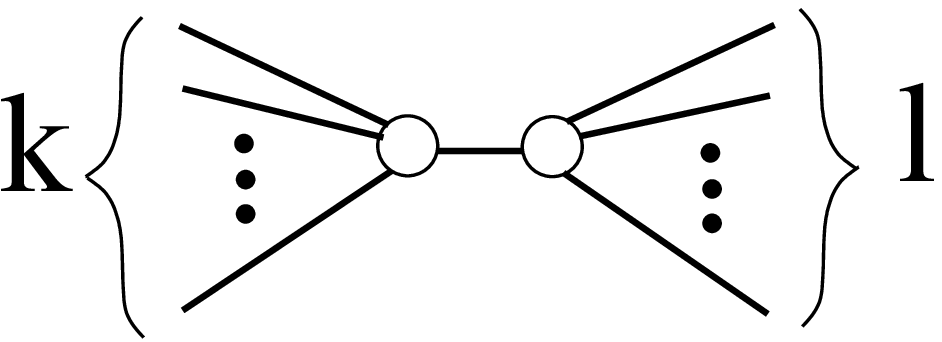,height=.3in}} & \to &
p\ \parbox{.6in}{\psfig{figure=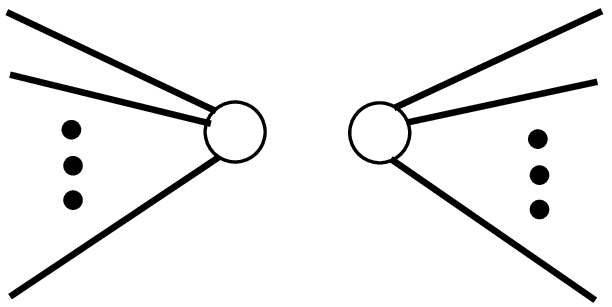,height=.3in}} + q\ 
\parbox{.6in}{\psfig{figure=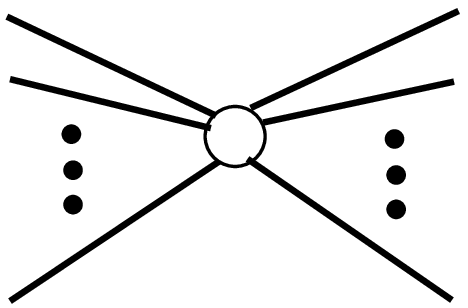,height=.3in}},\\

S_l:& \parbox{0.8in}{\psfig{figure=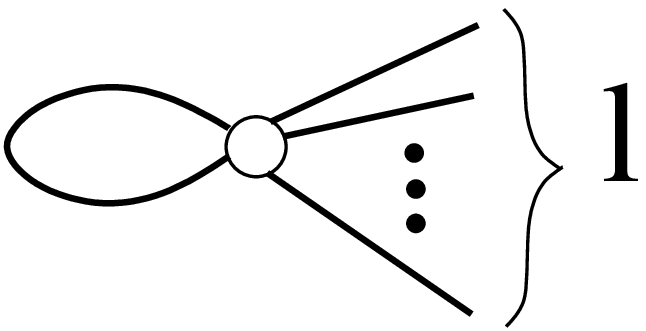,height=.3in}} & \to &
s\ \parbox{.3in}{\psfig{figure=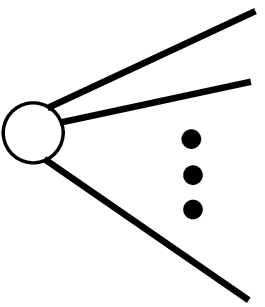,height=.3in}},\\

S_v: & \diag{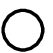}{.08in} & \to & v\emptyset,\\

S_{bw}: & \parbox{.5in}{\psfig{figure=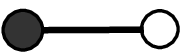,height=.08in}} & \to &
w_1\ \diag{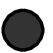}{.08in}\\ 

S_{bwb}: & \parbox{.5in}{\psfig{figure=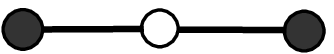,height=.08in}} & \to &
w_2\ \parbox{.5in}{\psfig{figure=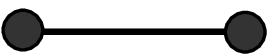,height=.08in}}\\
\end{array}$$
where external (respectively: internal) vertices are denoted by black 
(respectively: white) nodes and $k,l\in \Z_{\geq 0}$. 

\begin{theorem}\label{dichrom}
(a) The above reduction rules are terminal and confluent.\\
(b) Connected components of irreducible graphs are either isolated external
 vertices, \diag{blackdot}{.08in}, or 
\parbox{.3in}{\psfig{figure=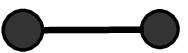,height=.08in}} or
at least $3$ external vertices connected to a single internal vertex.
\end{theorem}

The termination of the above rules is obvious. The proof of confluence 
is left to the reader.

\begin{corollary}
Irreducible graphs in $\G_n$ are in $1$-$1$ correspondence to partitions 
of $\{1,...,n\}$.
\end{corollary}

The only irreducible graph with no external 
vertices is $\emptyset$. Consequently, 
$R\G_0/\R(S_{kl}, S_l, k,l\geq 0, S_0)$ is a 
cyclic $R$-module generated by $\emptyset$. 
The projection $\G_0\to R\G_0/\R(S_{kl}, k,l\geq 0, S_0)\simeq R$
followed by substitution $R\stackrel{p\to 1, s\to 1+q}{\longrightarrow} 
\Z[q,v]$ is the dichromatic polynomial of graphs.
A generalization of dichromatic polynomial to ribbon graphs is considered in
\cite{BR}. It can be defined by reductions rules similar to those above as 
well.

\subsection{Partition Category}
Let $R$ be a ring with a specified $\delta^{\pm 1}\in R$.
A version of the ``dichromatic'' reduction rules appears in the context
of the partition category. There are two types of that category: symmetric 
and planar one. In each of them, objects are non-negative integers.

In the symmetric partition category, the morphisms $[n]\to [m]$
are $R$-linear combinations
of abstract graphs with $n+m$ external vertices divided into an
{\em input set} and an {\em output set} of $n$ and $m$ vertices, respectively. 
Vertices in each of these sets are numbered. The internal vertices 
have valency at most $3.$ 
Furthermore, the graphs are subject to the following rules:

$$T_1:\ \diag{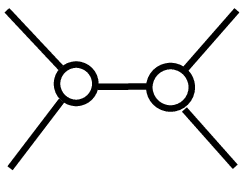}{.25in}\quad =\ \diag{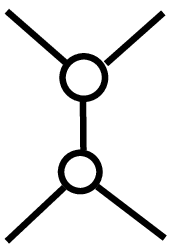}{.35in},\quad
T_2:\ \diag{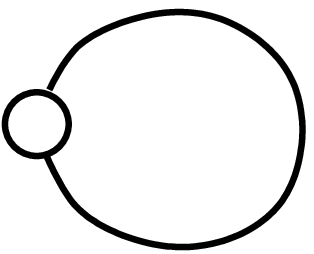}{.3in}\ =\ \delta,\quad
T_3: \parbox{.7in}{\psfig{figure=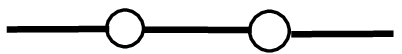,height=.08in}} =
\parbox{.5in}{\psfig{figure=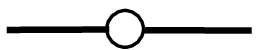,height=.08in}},$$
$$T_4:\ \parbox{.7in}{\psfig{figure=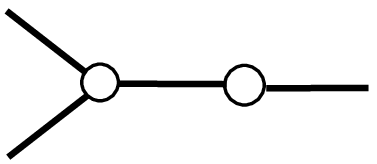,height=.3in}} =
\parbox{.5in}{\psfig{figure=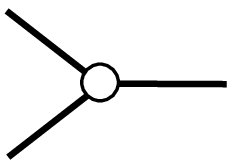,height=.3in}},\quad
T_5:\ \parbox{.4in}{\psfig{figure=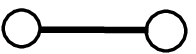,height=0.08in}}=\ 1,\ 
T_6:\ \parbox{.4in}{\psfig{figure=dichrom7,height=.08in}} =\quad
\diag{blackdot}{.08in},$$
$$T_7:\ \diag{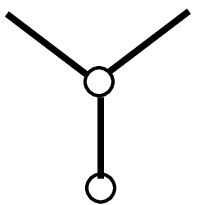}{.4in}\ =\ \diag{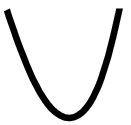}{.25in},\ 
T_8:\ \parbox{.25in}{\psfig{figure=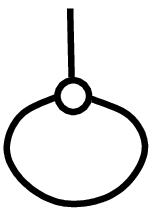,height=0.4in}}\ =
\ \delta \parbox{.2in}{\psfig{figure=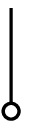,height=0.4in}}.$$

As before, empty/full dots represent internal/external vertices.
The composition of morphisms $\Gamma_1\in Mor([n],[m]),$
$\Gamma_2\in Mor([m],[k])$ is defined by identifying  
output vertices of $\Gamma_1$ with corresponding input vertices of $\Gamma_2.$
(The correspondence between these vertices is established by their numbering.)
These vertices become internal in $\Gamma_2\circ \Gamma_1.$
The symmetric partition category is monoidal, with the tensor product of 
morphisms $\Gamma_1\in Mor([n_1],[m_1]),$ $\Gamma_2\in Mor([n_2],[m_2])$ 
given by taking their disjoint union and shifting the numbering of the input 
and output vertices of $\Gamma_2$ by $n_1$ and $m_1$ respectively.

In the planar partition category, one thinks of objects as
sets $[n]=\{\frac{1}{n+1},...,\frac{n}{n+1}\}\subset [0,1].$
Morphisms $[n]\to [m]$ are $R$-linear (manifold) graphs in $[0,1]\times [0,1]$ 
without $S^1$'s, whose external vertices are 
$[n]\times \{0\}\cup [m]\times \{1\}.$
These graphs are considered up to relations $T_1,...,T_8.$
We denote the space of morphisms in the planar partition category by
$Mor^p([n],[m])$.
Compositions (respectively: tensor products) of morphisms are given by 
vertical (respectively: horizontal) stacking of graphs.

The importance of the planar (respectively: symmetric) partition category
stems from the fact that the object $[1]$ is the universal Frobenius 
algebra (respectively: commutative Frobenius 
algebra) in a monoidal (respectively symmetric monoidal) category.
Consequently, monoidal (respectively: symmetric monoidal) functors from the
planar (respectively: symmetric) 
partition category into the category of $R$-modules
are in $1$-$1$ correspondence with Frobenius algebras (respectively:
commutative Frobenius algebras) over $R$. Furthermore, the algebra of 
morphisms $[n] \to [n]$ called the {\em partition algebra} appears in the 
theory of Potts models, \cite{Jo,Ma}.

Reduction relations going from left to right sides of equations $T_1$-$T_8$
are locally confluent but not terminal since $T_1$ is invertible.
(Nonetheless, we are going to prove that they are confluent.)
The intuitive way of making these rules terminal is by allowing
$4$-valent vertices and by adding a reduction rule 
$$\diag{parcath}{.25in}\quad \to \diag{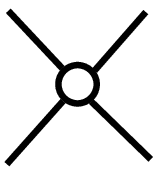}{.25in}.$$
Now to make the new rules terminal and locally confluent one needs 
to add an additional relation involving $5$-valent vertices, and then
one involving $6$-valent vertices, etc. Finally, one arrives at:

\begin{theorem}\label{part_cat_thm}
(1) The natural embedding of graphs of the symmetric partition category
into (abstract) graphs factors to an isomorphism of $R$-modules
$$\phi: Mor([n],[m]) \to \G_{n+m}/\R({\bar S}_{k,l},{\bar S}_l,{\bar S}_v,
{\bar S}_{bw},{\bar S}_{bwb}),$$
where
$${\bar S}_{k,l}:\ \parbox{0.8in}{\psfig{figure=dichrom1,height=.3in}}\ \to \ 
\parbox{.6in}{\psfig{figure=dichrom3,height=.3in}},\quad
{\bar S}_l:\ \parbox{0.8in}{\psfig{figure=dichrom9,height=.3in}}\ \to \ 
\delta\ \parbox{.3in}{\psfig{figure=dichrom10,height=.3in}},$$
$${\bar S}_v:\ \diag{whitedot}{.08in}\ \to\ \emptyset,\quad
{\bar S}_{bw}:\ \parbox{.4in}{\psfig{figure=dichrom7,height=.08in}}\ \to
\diag{blackdot}{.08in},\quad
{\bar S}_{bwb}:\ \parbox{.5in}{\psfig{figure=dichrom4,height=.08in}}\ \to \ 
\parbox{.5in}{\psfig{figure=dichrom6,height=.08in}}.$$
(Note that these are the dichromatic reduction rules for 
$p=0, q=v=w_1=w_2=1, s=\delta$).\\
(2) Similarly, the natural embedding of graphs of planar 
partition category into planar graphs in $D^2$ with $n+m$ external vertices
factors to an isomorphism 
$$\phi^p: Mor^p([n],[m]) \to \G_{n+m}(D^2)/
\R(\bar S_{k,l},\bar S_l,\bar S_v,\bar S_{bw},\bar S_{bwb}).$$
\end{theorem}

\begin{proof}[Sketch of Proof]
(1) The irreducible graphs in $\G_{n+m}$ listed in 
Theorem \ref{dichrom} span $\G_{n+m}/\R({\bar S}_{k,l},{\bar S}_l,{\bar S}_v,
{\bar S}_{bw},{\bar S}_{bwb}).$ Since all of these graphs are values of
$\phi,$ it is an epimorphism.
To prove that $\phi$ is $1$-$1$ observe that connected components of 
every graph in $Mor([n],[m])$ determine a partition of $\{1,...,n+m\}.$
Let $Mor^\tau([n],[m])$ be the subspace of $Mor([n],[m])$
spanned by graphs associated with the partition $\tau.$
Since $T_1,...,T_8$ preserve partitions,
$$Mor([n],[m])=\bigoplus_{\text{partitions }\tau} Mor^\tau([n],[m])$$
and, similarly, $\G_{n+m}/\R(\bar S_{k,l}, \bar S_l, \bar S_v, \bar S_{bw},
\bar S_{bwb})$ decomposes into subspaces indexed
by partitions, which by Theorem \ref{dichrom} are $1$-dimensional.
Since $\phi$ preserves partition classes, it is enough to prove that
$Mor^\tau([n],[m])\simeq R$ as an $R$-module.
This follows from the following:

\begin{lemma}
Any two connected graphs with internal vertices of valency $\leq 3$ 
are equivalent via relations $T_1,...,T_8$.
\end{lemma}

\noindent {\it Proof:} $T_4$ and $T_1$ allow to ``slide'' edges past 
$2$-valent and $3$-valent vertices. Therefore,
all cycles in a graph can be transformed into loops, which can be 
eliminated by $T_2$ and $T_8$. Furthermore, all internal $1$-valent 
vertices can be removed by $T_5, T_6$ and $T_7$. Consequently, each 
connected graph is equivalent to a tree with no internal $1$-valent vertices. 
All such trees are related by $T_1$ moves.
There is a geometric way to see that. First, any such tree embeds into $D^2$
and its dual corresponds to a
division of a $(n+m)$-gon into triangles by non-intersecting diagonals.
Here is an example of a planar tree and the corresponding dual 
triangulation of $6$-gon:
$$\diag{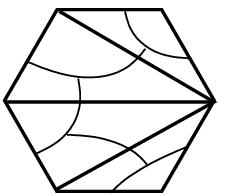}{.6in}$$

Any two such triangulations are related by the move:

$$\diag{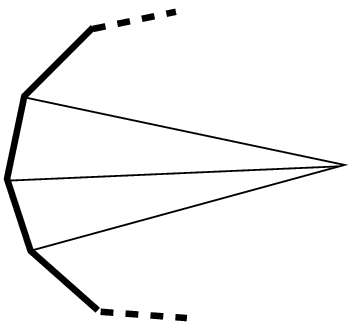}{.7in}\quad \leftrightarrow\quad \diag{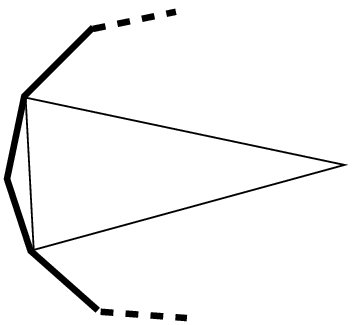}{.7in}$$
which is dual to $T_1$.
\end{proof}

The proof of Theorem \ref{part_cat_thm}(2) is analogous.

%

%
\section{Application to knots}\label{s_knots}
%

Now, we turn to spaces of dimension $3$, which are the most difficult
to deal with in the context of graph embeddings.
Graphs in a $3$-dimensional manifold $M$ are the easiest to analyze if 
$M$ is an $I$-bundle over a surface $F$, $I=[-1,1],$ since then each 
graph is represented by its diagram in $F$ and such representations are
unique up to Reidemeister moves:

$$RI:\ \diag{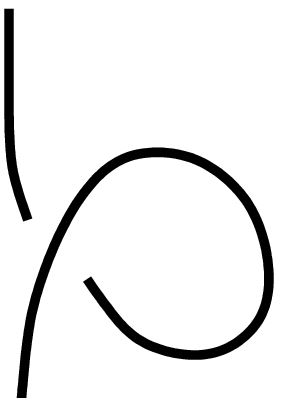}{.4in} \leftrightarrow 
\parbox{.2in}{\psfig{figure=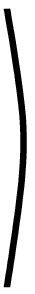,height=.4in}},
\quad RII:\ \diag{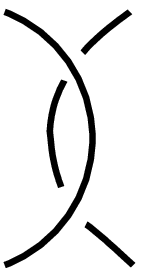}{.3in}\leftrightarrow  \diag{smootha}{.3in},
\quad \text{and}\quad 
RIII:\ \diag{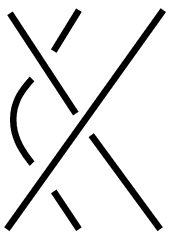}{.3in}\leftrightarrow \diag{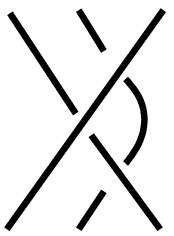}{.3in},$$
and the moves 
$$V: \quad 
\parbox{.7in}{\psfig{figure=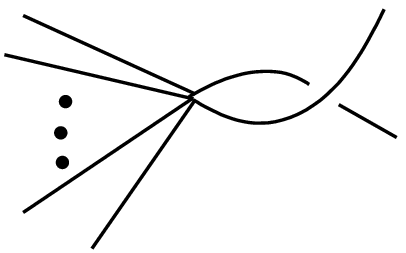,height=0.4in}}
\ \leftrightarrow \ 
\parbox{.7in}{\psfig{figure=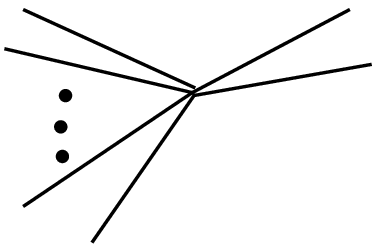,height=0.4in}}
\ \leftrightarrow \ 
\parbox{.7in}{\psfig{figure=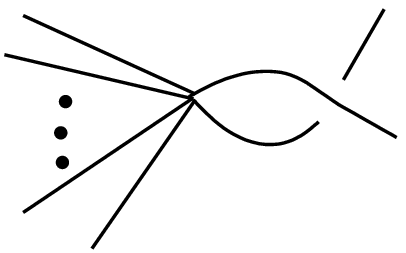,height=0.4in}}$$
$$V_{k,l}:\quad \diag{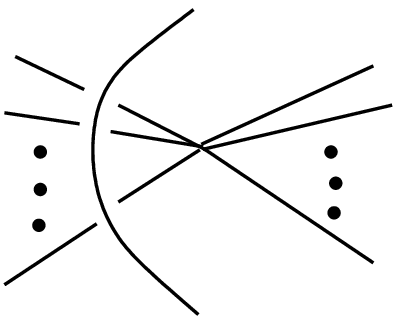}{.5in}\ \quad 
\leftrightarrow \quad \diag{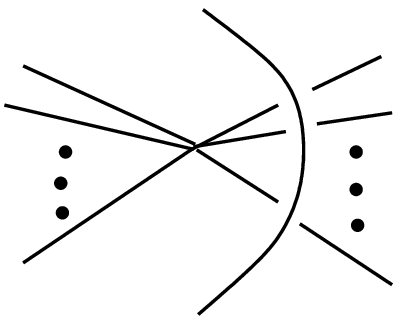}{.5in}$$ where there are 
$k$ edges on the right and
$l$ on the left, for all $k,l\geq 0$.

The main problem in knot theory is deciding whether any two
link diagrams represent isotopic links. This reduces now to question
whether these link diagrams are equal in 
$\Z\ld(F)/\R(RI,RII,RIII),$ 
where $\ld(F)$ is the set of all link diagrams in $F$. 
(In classical knot theory $F=D^2$, but other 
surfaces are of interest for us as well.)
Notice that this is a version of the
problem formulated in Introduction.

The rules $$r_1:\ \diag{ri}{.4in}\to
\parbox{.2in}{\psfig{figure=linev.eps,height=.4in}},\quad
r_2:\ \diag{rii}{.3in} \to \diag{smootha}{.3in}$$
although terminal are not confluent, since 
$$\diag{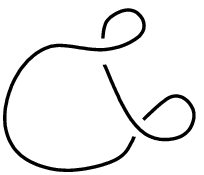}{.35in}\quad \text{and}\ \diag{circle}{.3in}$$
are both irreducible but not equivalent to each other.
However, it is easy to show by the method of Section \ref{ss_overlaps}
that $r_1,r_2$ together with
$$r_1': \ \diag{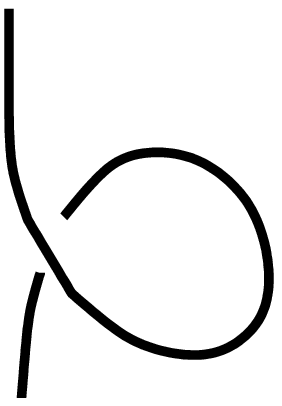}{.4in}\to
\parbox{.2in}{\psfig{figure=linev.eps,height=.4in}}$$
are terminal and confluent.
Unfortunately, the rule $r_3: \diag{riii1}{.3in} \to \diag{riii2}{.3in}$
is not terminal since it is its own inverse. Furthermore, 
$r_1,r_1',r_2,r_3$ are not confluent.

Nonetheless, inspired by the notion of confluence 
we ask the following question:

\begin{question} For any given surface $F,$ is there a set $\G$ of 
graphs in $F$ containing $\ld(F)$ and a 
finite set of terminal and confluent reduction rules $S_1,...,S_d,$
with coefficients in a ring $R$ 
such that the inclusion $\ld(F)\hookrightarrow \G$ induces a monomorphism
$$\phi: R\ld(F)/\R(RI,RII,RIII)\hookrightarrow R\G/\R(S_1,...,S_d).$$
\end{question}

A positive answer to this question would provide an obvious algorithm for
distinguishing non-isotopic links in $I$-bundles over $F$. We conjecture 
that there does not exist a set of terminal and confluent reduction
rules with these properties. 

Confluence theory provides an immediate proof of the following statement.

\begin{theorem}
Let $F$ be any surface and $B\subset \p F$ a finite set.\\
(1) unoriented link diagrams ($A_1$-webs) in $F$ are invariant in 
${\mathbb A}_1(F,B,R)$ under moves 
$RII,$ $RIII,$ and the first balanced Reidemeister move: 
$$RIb: \diag{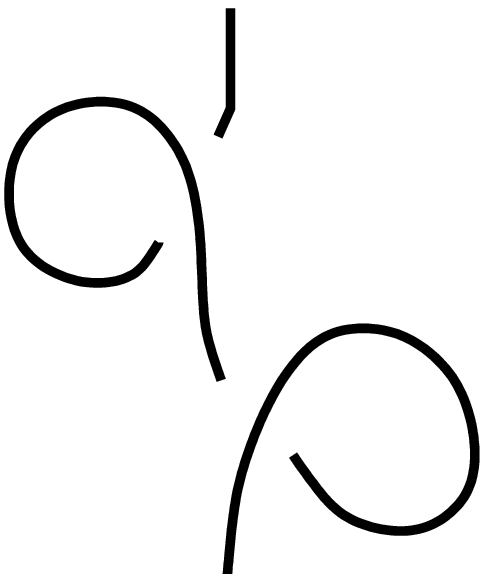}{.4in}\to 
\parbox{.2in}{\psfig{figure=linev.eps,height=.4in}}.$$
(2) $A_2$-webs are invariant in ${\mathbb A}_2(F,B,R)$ under
oriented $2$nd and $3$rd Reidemeister moves as well as 
oriented $1$st balanced Reidemeister move and 
under moving an arc over a vertex.\\
(3) Similarly, $B_2$-webs and $G_2$-webs in $F$ are invariant under 
all 2nd and 3rd Reidemeister moves (involving both single and double 
lines) and under moving an arc over a vertex.
\end{theorem}

To prove invariance under any of the above relations it is enough to 
check that reduction rules applied to both sides of that relation yield 
identical linear diagrams.

\begin{corollary} For orientable $F$, 
$\A_1(F,\emptyset,R), \B_2(F,\emptyset,R),\G_2(F,\emptyset,R),$ 
provide invariants of framed unoriented links and 
$\A_2(F,\emptyset,R)$ provides an invariant of framed oriented links 
under isotopy in $F\times I$.
\end{corollary} 

For $F=\Rb^2,$ $B=\emptyset,$ and $R=\Z[A^{\pm 1}],$ the module 
$\A_1(F,B,R)$ is free on one generator, $\emptyset,$ 
and $[L]\in \A_1(F,B,R)=R$ is the Kauffman bracket of $L,$ \cite{Ka}.
For other oriented surfaces, $\A_1(F,\emptyset,R)$ is isomorphic to the 
Kauffman bracket skein module of $F\times I,$ \cite{HP2, P-intro, 
P-fundamentals, PS}. (For more on Kauffman bracket skein modules see 
\cite{Bu-knot,Bu-char, BFKu, BFKym, BHMV, BP, FG, FGL, FK, GS1, GS2, GH,
HP3, HP4, Le,Sa1, Sa2, S-TQFT, S-4th, Tu}.)
Consequently, Corollary \ref{A_1basis} immediately implies
the result of Przytycki, \cite[Theorem 3.1]{P-fundamentals}:

\begin{theorem}\label{kb-surf1}
For any ring $R$ 
with $A^{\pm 1}\in R$ and any orientable surface $F,$ the Kauffman 
bracket skein module of $F\times [0,1]$ is a free $R$-module with 
basis composed of links whose diagrams in $F$ have no crossings and no 
contractible components.
\end{theorem}

There are versions of this theorem for orientable $I$-bundles over 
non-orientable surfaces and for $B\ne \emptyset$. They can be easily proved
by the method of confluence as well.

The module $\A_2(F,\emptyset,R)$ is isomorphic to the $SU_3$-skein module of
$F\times I$ introduced in \cite{FZ, S-SUn}. (See also \cite{OY}.)
Consequently, Corollary \ref{A_2basis} implies

\begin{theorem}
The $SU_3$-skein module of $F\times I,$ $\S_3(F\times I,R)$ (in notation of
\cite{S-SUn}) is a free $R$-module with a basis given by all $A_2$-webs in $F$
with no $0$-gons, no true bigons, and no true $4$-gons.
\end{theorem}

\begin{problem}
A large number of skein modules is considered in the literature,
other than those mentioned above,
\cite{AT1,AT2,GZ,HM,HP1,Kai1, Kai2, Kai3, KL, Li, P-handlebody,P-PAN, P-vasgus,
  P-qanalog,P-deformations, P-homotopy,P-cubic,PT, Zh,ZL}. Can the method of 
confluence be 
applied to determine canonical bases of these modules for $F\times I$?
\end{problem}

The applications of confluence theory to knot theory discussed so far 
apply to links in $I$-bundles over surfaces. Unfortunately, reduction 
rules for links in arbitrary $3$-manifolds are more difficult to handle.
This is illustrated by the Kauffman bracket skein relations:

Let $\mathcal L(M)$ be the set of all framed unoriented links in an
orientable $3$-manifold $M$. Let $R$ be a ring with a distinguished element 
$A^{\pm 1}\in R$ and let
$$\begin{array}{llll}
S_1: &  \diag{crossa}{.3in} & \to & A\diag{smootha}{.3in}+A^{-1}
\diag{smoothb}{.3in}\\
S_2: & \diag{circle}{.3in} & \to & (A^2+A^{-2})\emptyset
\end{array}$$
be reduction rules taking place in $D^3$.
The $R$-module
$R{\mathcal L}(M)/\R(S_1,S_2)$ is called {\em the Kauffman bracket skein module
of $M$}. We leave the proof of the following to the reader:

\begin{proposition}
$S_1,S_2$ are confluent but not terminal.
\end{proposition}

%

\end{document}